\documentclass[12pt, a4paper]{article}
\usepackage[latin1]{inputenc}
\usepackage{amsmath}
\usepackage[english]{babel}
\usepackage{indentfirst}
\usepackage{amstext}
\usepackage{amsfonts}
\usepackage{textcomp}
\usepackage{color}
\usepackage{amssymb}
\usepackage[dvips]{graphicx}
\usepackage{setspace}
\usepackage[all]{xy}

\newcommand {\demo}{\hskip -0.6cm{\bf Proof.  }}
\newcommand {\fim}{\hfill{$\square$}\vskip 1pc}

\newcommand {\R}{\mathbb{R}}
\newcommand {\N}{\mathbb{N}}

\newcommand{\supp}{\text{supp}}

\newcommand{\id}{\mathrm{id}}

\newtheorem{teorema}{Theorem}[section]
\newtheorem{lema}[teorema]{Lemma}
\newtheorem{corolario}[teorema]{Corollary}
\newtheorem{definicao}[teorema]{Definition}
\newtheorem{proposicao}[teorema]{Proposition}
\newtheorem{exemplo}[teorema]{Example}
\newtheorem{rmk}[teorema]{Remark}

\begin{document}

\onehalfspace

\title{Branching Systems and General Cuntz-Krieger Uniqueness Theorem for Ultragraph C*-algebras}
\author{Daniel Gon\c{c}alves\footnote{Partially supported by CNPq.}, Hui Li\footnote{Supported by Research Center for Operator Algebras of East China Normal University. The corresponding author.} \ and Danilo Royer}
\date{26 Jan 2016}
\maketitle

AMS 2000 MSC: 47L99, 37A55


\begin{abstract}
We give a notion of branching systems on ultragraphs. From this we build concrete representations of ultragraph C*-algebras on the bounded linear operators of Hilbert spaces. To each branching system of an ultragraph we describe the associated Perron-Frobenius operator in terms of the induced representation. We show that every permutative representation of an ultragraph C*-algebra is unitary equivalent to a representation arising from a branching system. We give a sufficient condition on ultragraphs such that a large class of representations of the C*-algebras of these ultragraphs is permutative. To give a sufficient condition on branching systems so that their induced representations are faithful we generalize Szyma{\'n}ski's version of the Cuntz-Krieger uniqueness theorem for ultragraph C*-algebras.
\end{abstract}

\section{Introduction}

Ultragraphs are combinatorial objects that generalize directed graphs. Roughly speaking an ultragraph is a graph where the image of the range map does not belong to the set of vertices, but instead to the power set on vertices. The concept was introduced by Mark Tomforde in \cite{Tomforde:JOT03} with an eye towards C*-algebra applications. In particular, Tomforde showed in \cite{Tomforde:JOT03} how to associate a C*-algebra to an ultragraph and proved that there exist ultragraph C*-algebras that are neither Exel-Laca algebras nor graph C*-algebras. So the study of ultragraph C*-algebras is of interest, but furthermore, ultragraph C*-algebras were key in answering the long-standing question of whether an Exel-Laca algebra is Morita equivalent to a graph algebra (see \cite{KMST}).

Recently the scope of ultragraphs has surpassed the realm of C*-algebras, reaching applications to symbolic dynamics. In particular, ultragraphs are fundamental in characterizing one-sided shift spaces over infinite alphabets. Furthermore, questions regarding the dynamics of one-sided shift spaces over infinite alphabets were answered using ultragraphs and their C*-algebras (see \cite{GRUltragraph}).  

The above evidence leads us to believe that there are still many applications of ultragraphs to be found. Indeed, among the results we present in this paper, we will show a connection (via branching systems) between ultragraphs and the Perron-Frobenius operator from the ergodic theory (see Section~\ref{PFoperator}), therefore generalizing results previously obtained for graph algebras and Cuntz-Krieger C*-algebras (see \cite{MR2903145, GR}). 

Branching systems are not just important as a way to connect ultragraphs to the ergodic theory. They have appeared in fields as random walks, symbolic dynamics, wavelet theory and are strongly connected to the representation theory of combinatorial algebras.  In particular, Bratteli and Jorgensen have initiated the study of wavelets and representations of the Cuntz algebra via branching systems in \cite{MR1465320, BJ}. After this, many results relating branching systems and representations of generalized Cuntz algebras were obtained, see for example \cite{FGKP, GLR, GR4, GR3, GR, MR2848777}. 
It is our goal in this paper to generalize many of the results obtained in the just mentioned manuscripts to ultragraphs.

In particular, to obtain ultragraph versions of the results presented in \cite{GLR} we must first develop a generalized Cuntz-Krieger uniqueness theorem for ultragraph C*-algebras. Actually there has been a lot of activity recently on generalized Cuntz-Krieger theorems for combinatorial algebras, see for example \cite{BNS,BNSSW, Szyma'nski:IJM02} and so the proof of a generalized Cuntz-Krieger uniqueness theorem for ultragraph C*-algebras is of independent interest.

We break the paper as follows. In Section~2 we present a preliminary on ultragraphs. In Sections~3, 4 and 5 we generalize the results in \cite{MR2903145} to ultragraphs. In particular we introduce branching systems of ultragraphs in Section~3; In Section~4 we show how to build a representation of the ultragraph C*-algebra from a branching system of an ultragraph; then in Section~5 we find the connection between branching systems, representations and the Perron-Frobenius operator. We generalize the results in \cite{MR2848777} in Section~6, where we introduce permutative representations, show that they are always unitarily equivalent to a representation arising from a branching system and give sufficient conditions on ultragraphs so that a large class of representations on these ultragraph C*-algebras are permutative. In Section~7 we make a pause in the theory of branching systems and prove the generalized Cuntz-Krieger uniqueness theorem for ultragraph C*-algebras. This section can be read independently from the other sections. To end the paper, in Section~8 we prove versions of the results in \cite{GLR} for ultragraphs, that is, we prove a converse of the Cuntz-Krieger uniqueness theorem for ultragraph C*-algebras and give a sufficient condition for a representation of an ultragraph C*-algebra, arising from a branching system, to be faithful.

\section{Preliminaries}

Throughout this paper, the notation $\mathbb{N}$ stands for the set of all the positive integers; all measure spaces are assumed to be $\sigma$-finite; all C*-algebras are assumed to be separable; and all representations of C*-algebras are assumed to be on separable infinite-dimensional Hilbert spaces.

The following lemma might be well-known but we are not able to find any reference to it.

\begin{lema}\label{imprimitivity bimodule}
Let $A$ be a C*-algebra, let $I,J$ be closed two-sided ideals of $A$, let $B$ be a C*-subalgebra of $I$, and let $X$ be a closed subspace of $I$. Suppose that $BX,XI \subset X; XX^* \subset B,X^*X \subset I$; and $\overline{\mathrm{span}}XX^*=B,\overline{\mathrm{span}}X^*X=I$. Furthermore, suppose that $B \cap J=0$. Then $I \cap J=0$.
\end{lema}

\demo
It is straightforward to show that $X$ is a $B$--$I$ imprimitivity bimodule (see \cite[Definition~3.1]{RaeburnWilliams:Moritaequivalenceand98}) with the natural operations on $X$. By \cite[Theorem~3.22]{RaeburnWilliams:Moritaequivalenceand98}, there exists a lattice isomorphism $\Phi$ from the set of closed two-sided ideals of $B$ onto the set of closed two-sided ideals of $I$. A straightforward calculation gives
\[
\Phi(B \cap J)=\overline{\mathrm{span}}\{x^*y:x,y \in X; \text{ and for } z \in X, \text{ we have } yz^*\in J \cap B \}.
\]
We claim that $\Phi(B \cap J)=I \cap J$. Since $X \subset I, \Phi(B \cap J)\subset I$. Fix $x,y \in X$ satisfying that for $z \in X,yz^*\in J \cap B$. Notice that $x^*y(\sum_{i=1}^{n}z_i^*z_i') \in J$ for $z_1,z_1',\dots,z_n,z_n' \in X$. Since $\overline{\mathrm{span}}X^*X=I$, we deduce that $x^*y \in J$. So $\Phi(B \cap J) \subset I \cap J$. Conversely, fix $\sum_{i=1}^{n}x_i^*x_i' \in \mathrm{span}X^*X \cap J$. Take an approximate identity $(e_\alpha) \subset I \cap J$. Then $x_i' e_\alpha z^* \in J \cap B$ for all $\alpha$, for all $z \in X$ , and $\sum_{i=1}^{n}x_i^*x_i' e_\alpha \to \sum_{i=1}^{n}x_i^*x_i'$. So $\sum_{i=1}^{n}x_i^*x_i' \in \Phi(B \cap J)$. Since $\overline{\mathrm{span}}X^*X=I, I \cap J\subset \Phi(B \cap J)$. Hence $\Phi(B \cap J)=I \cap J$. Therefore $I \cap J=0$ because $B \cap J=0$. \fim

\begin{rmk}
$\Phi$ as in Lemma~\ref{imprimitivity bimodule} is indeed the inverse map of the Rieffel correspondence given in \cite[Proposition~3.24]{RaeburnWilliams:Moritaequivalenceand98}. 
\end{rmk}

In the rest of this section we give a brief introduction to ultragraphs and ultragraph C*-algebras, as defined by Tomforde in \cite{Tomforde:JOT03}.

\begin{definicao}[{\cite[Definition~2.1]{Tomforde:JOT03}}]\label{def of ultragraph}
An \emph{ultragraph} is a quadruple $\mathcal{G}=(G^0, \mathcal{G}^1, r,s)$ consisting of two countable sets $G^0, \mathcal{G}^1$, a map $s:\mathcal{G}^1 \to G^0$ and a map $r:\mathcal{G}^1 \to P(G^0)\setminus \{\emptyset\}$, where $P(G^0)$ stands for the power set of $G^0$.
\end{definicao}

\begin{exemplo}\label{exeUltragrafo} Below we show the picture of the ultragraph where $G^0 = \{v_1, \ldots v_{10}\}$, $\mathcal{G}^1 = \{e_1, \ldots, e_5\}$ and the source and range maps are defined as follows: $s(e_1) = v_1$, $s_(e_2)=v_6$, $s(e_3)=v_2$, $s(e_4)= v_6$ and $s(e_5)=v_{10}$; $r(e_1)=\{v_2,v_3 \}$, $r(e_2)=\{v_3,v_4,v_5 \}$, $r(e_3)=\{v_7 \}$, $r(e_4)=\{ v_{10}\}$ and $r(e_5)=\{ v_8, v_9 \}$. 

\vspace{1cm}

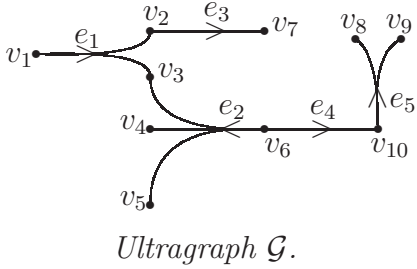
\centerline{
\setlength{\unitlength}{1cm}
\begin{picture}(11,0)
\put(-0.5,0){\circle*{0.1}}
\put(-0.9,-0.1){$v_1$}
\qbezier(-0.5,0)(1,0.05)(1,-0.3)
\qbezier(-0.5,0)(1,-0.05)(1,0.3)
\put(1,0.3){\circle*{0.1}}
\put(0.9,0.45){$v_2$}
\put(1,-0.3){\circle*{0.1}}
\put(1.1,-0.3){$v_3$}
\put(0,-0.1){$>$}
\put(0,0.15){$e_1$}
\put(2.5,-1){\circle*{0.1}}
\put(2.5,-1.3){$v_6$}
\qbezier(2.5,-1)(1,-1.05)(1,-0.3)
\qbezier(2.5,-1)(1,-1)(1,-1)
\qbezier(2.5,-1)(1,-0.9)(1,-2)
\put(1,-1){\circle*{0.1}}
\put(0.6,-1){$v_4$}
\put(1,-2){\circle*{0.1}}
\put(0.6,-2){$v_5$}
\put(1.9,-1.1){$<$}
\put(1.9,-0.8){$e_2$}
%\put(2.5,0){\circle*{0.1}}
%\put(2.1,0){$v_7$}
\put(2.5,0.3){\circle*{0.1}}
\put(2.6,0.3){$v_7$}
\put(1,0.3){\line(1,0){1.5}}
\put(1.7,0.2){$>$}
\put(3.1,-0.8){$e_4$}
\put(4,-1){\circle*{0.1}}
\put(3.9,-1.3){$v_{10}$}
\put(2.5,-1){\line(1,0){1.5}}
\put(3.1,-1.1){$>$}
\put(1.7,0.5){$e_3$}
\qbezier(4,-1)(3.9,0)(4.3,0.2)
\qbezier(4,-1)(4,0)(3.7,0.2)
\put(3.85,-0.7){$\wedge$}
\put(4.15,-0.7){$e_5$}
\put(3.7,0.2){\circle*{0.1}}
\put(3.5,0.35){$v_8$}
\put(4.3,0.2){\circle*{0.1}}
\put(4.1,0.35){$v_9$}
\put(0.5,-2.7){Ultragraph $\mathcal{G}$. }
\end{picture}}

\end{exemplo}
\vspace{2.5cm}

\begin{rmk} Notice how in the drawing above we split the edges to represent their range.
\end{rmk}

\begin{definicao}[{\cite[Page~349]{Tomforde:JOT03}}]\label{defofultragraph}
Let $\mathcal{G}$ be an ultragraph. Define $\mathcal{G}^0$ to be the smallest subset of $P(G^0)$ that contains $\{v\}$ for all $v\in G^0$, contains $r(e)$ for all $e\in \mathcal{G}^1$, and is closed under finite unions and nonempty finite intersections.
\end{definicao}

The set $\mathcal{G}^0$ can be characterized in the following way.

\begin{lema}[{\cite[Lemma~2.12]{Tomforde:JOT03}}]\label{concrete des of G^0}
Let $\mathcal{G}$ be an ultragraph. Then
\begin{align*}
\mathcal{G}^0=\Big\{\Big(\bigcap_{e\in X_1}r(e)\Big) \cup\dots\cup\Big(\bigcap_{e\in X_n}r(e)\Big) \cup F&:X_i \text{'s are finite subsets of } \mathcal{G}^1,
\\&F \text{ is a finite subset of } G^0\Big\}.
\end{align*}
\end{lema}

\begin{definicao}[{\cite[Definition~2.7, Theorem~2.11]{Tomforde:JOT03}}]\label{defofultragraphalgebra}
Let $\mathcal{G}$ be an ultragraph. The \emph{ultragraph algebra} $C^*(\mathcal{G})$ is the universal C*-algebra generated by a family of partial isometries with orthogonal ranges $\{s_e:e\in \mathcal{G}^1\}$ and a family of projections $\{p_A:A\in \mathcal{G}^0\}$ satisfying
\begin{enumerate}
\item\label{p_Ap_B=p_{A cap B}} $p_\emptyset=0, p_Ap_B=p_{A\cap B}, p_{A\cup B}=p_A+p_B-p_{A\cap B}$, for all $A,B\in \mathcal{G}^0$;
\item\label{s_e^*s_e=p_{r(e)}}$s_e^*s_e=p_{r(e)}$, for all $e\in \mathcal{G}^1$;
\item $s_es_e^*\leq p_{s(e)}$ for all $e\in \mathcal{G}^1$; and
\item\label{CK-condition} $p_v=\sum\limits_{s(e)=v}s_es_e^*$ whenever $0<\vert s^{-1}(v)\vert< \infty$.
\end{enumerate}
Moreover, any family of partial isometries with orthogonal ranges $\{S_e:e\in \mathcal{G}^1\}$ and any family of projections $\{P_A:A\in \mathcal{G}^0\}$ in any C*-algebra $B$ satisfying Conditions~(\ref{p_Ap_B=p_{A cap B}})--(\ref{CK-condition}) is called a \emph{Cuntz-Krieger $\mathcal{G}$-family}.
\end{definicao}

Notice that $p_v \geq \sum_{e \in S}s_es_e^*$ for any nonempty finite subset $S \subset s^{-1}(v)$. It follows from \cite[Theorem~2.11]{Tomforde:JOT03} that each $s_e$ and $p_A$, with $A \neq \emptyset$, are nonzero.

\begin{definicao}[{\cite[Page~350]{Tomforde:JOT03}}]
Let $\mathcal{G}$ be an ultragraph. For $\alpha \in \mathcal{G}^0$, define $s(\alpha)=r(\alpha):=\alpha$, and define $\vert\alpha\vert:=0$. For $n \geq 1$, define $\mathcal{G}^n:=\{\alpha=(\alpha_i)_{i=1}^{n} \in \prod_{i=1}^{n}\mathcal{G}^1: s(\alpha_{i+1}) \in r(\alpha_i), i=1,\dots,n-1\}$, and for $\alpha \in \mathcal{G}^n$ define $s(\alpha):=s(\alpha_1), r(\alpha):=r(\alpha_n),\vert\alpha\vert:=n$. Define $\mathcal{G}^*:=\amalg_{n=0}^{\infty}\mathcal{G}^n$. For $\alpha \in \mathcal{G}^*$, define
\begin{align*}
s_\alpha:= \begin{cases}
    s_{\alpha_1}\dots s_{\alpha_{\vert\alpha\vert}}&\text{ if $\vert\alpha\vert >0$} \\
    p_{\alpha}&\text{ if $\vert\alpha\vert=0$.} \\
\end{cases}
\end{align*}
\end{definicao}

\begin{definicao}[{\cite[Definition~3.4]{Tomforde:JOT03}}]\label{defofConditionLforultragraph}
Let $\mathcal{G}$ be an ultragraph. An element $\alpha\in\mathcal{G}^*\setminus\mathcal{G}^0$ is called a \emph{cycle} if $s(\alpha_1) \in r(\alpha_n)$. A cycle $\alpha$ is said to be \emph{simple} if $\alpha_i \neq \alpha_j$ for all $i\neq j$. A cycle $\alpha$ is said to have \emph{exits} if one of the following happens:
\begin{enumerate}
\item there exists $1 \leq i \leq n-1$ such that $s^{-1}(r(\alpha_i)) \neq \alpha_{i+1}$;
\item $s^{-1}(r(\alpha_n)) \neq \alpha_{1}$;
\item there exist $1 \leq i \leq n$ and $v \in r(\alpha_i)$ such that $s^{-1}(v) =\emptyset$.
\end{enumerate}
Moreover, the ultragraph $\mathcal{G}$ is said to satisfy \emph{Condition~(L)} if every cycle has exits.
\end{definicao}

It is straightforward to see that a cycle $\alpha=(\alpha_1,\dots,\alpha_n)$ has no exits if $\vert r(\alpha_i) \vert=1, s^{-1}(s(\alpha_i))=\{\alpha_i\}$ for all $i=1,\dots,n$.

\begin{definicao}[{\cite[Definitions~3.1, 3.2]{Tomforde:IUMJ03}}]
Let $\mathcal{G}$ be an ultragraph. A subset $\mathcal{H} \subset \mathcal{G}^0$ is said to be \emph{hereditary} if
\begin{enumerate}
\item for $A,B \in \mathcal{H}$, we have $A \cup B \in \mathcal{H}$;
\item for $A \in \mathcal{H},B \in \mathcal{G}^0$, if $B \subset A$ then $B \in \mathcal{H}$; and
\item for $e \in \mathcal{G}^1$, if $\{s(e)\} \in \mathcal{H}$ then $r(e) \in \mathcal{H}$.
\end{enumerate}
A subset $\mathcal{S} \subset \mathcal{G}^0$ is said to be \emph{saturated} if for $v \in G^0, 0<\vert s^{-1}(v)\vert<\infty, r(s^{-1}(v)) \subset \mathcal{S} \implies \{v\} \in \mathcal{S}$.
\end{definicao}

\begin{lema}[{\cite[Lemma~3.12]{Tomforde:IUMJ03}}]\label{smallest here and sat}
Let $\mathcal{G}$ be an ultragraph and let $\mathcal{H}$ be a hereditary subset of $\mathcal{G}^0$. Define $\mathcal{H}_0:=\mathcal{H}$. For $n \geq 0$, define
\[
S_n:=\{v \in G^0:0<\vert s^{-1}(v)\vert <\infty, r(s^{-1}(v)) \subset \mathcal{H}_n \}; \text{ and define }
\]
\[
\mathcal{H}_{n+1}:=\{A\cup F:A \in \mathcal{H}_n,F \text{ is a finite subset of } S_n\}.
\]
Then $\bigcup_{n=0}^{\infty}\mathcal{H}_n$ is the smallest hereditary and saturated subset of $\mathcal{G}^0$ containing $\mathcal{H}$.
\end{lema}

\begin{lema}[{\cite[Lemma~3.5]{Tomforde:IUMJ03}}]\label{define ideal I_H}
Let $\mathcal{G}$ be an ultragraph and let $\mathcal{SH}$ be a hereditary and saturated subset of $\mathcal{G}^0$. Then the closed two-sided ideal $I(\mathcal{SH})$ of $C^*(\mathcal{G})$ generated by $\{p_A:A \in \mathcal{SH}\}$ is gauge-invariant and has the following form.
\[
I(\mathcal{SH})=\overline{\mathrm{span}}\{s_\alpha p_A s_\beta^*:\alpha,\beta\in\mathcal{G}^*,A \in\mathcal{SH}\}.
\]
\end{lema}

The following is the Cuntz-Krieger uniqueness theorem for ultragraph C*-algebras.

\begin{teorema}[{\cite[Theorem~6.7]{Tomforde:JOT03}}]\label{Cuntz-Krieger uni thm for ultragraph alg}
Let $\mathcal{G}$ be an ultragraph which satisfies Condition~(L). Then for a Cuntz-Krieger $\mathcal{G}$-family $\{P_A,S_e:A\in \mathcal{G}^0, e\in \mathcal{G}^1\}$ in a C*-algebra $B$ satisfying $P_A\neq 0$ whenever $A \neq \emptyset$, the ultragraph C*-algebra $C^*(\mathcal{G})$ is isomorphic with $C^*(P_A,S_e)$ via the map $p_A \mapsto P_A, s_e \mapsto S_e$.
\end{teorema}

%=============================================================================================================================
\section{Branching Systems of Ultragraphs}\label{BSU}
%=============================================================================================================================

In this section we introduce branching systems associated to ultragraphs and show that they always exist.

\begin{definicao}\label{branchsystem}
Let $\mathcal{G}$ be an ultragraph, $(X,\mu)$ be a measure space and let $\{R_e,D_A\}_{e\in \mathcal{G}^1,A\in \mathcal{G}^0}$ be a family of measurable subsets of $X$. Suppose that
\begin{enumerate}
\item\label{R_e cap R_f =emptyset if e neq f} $R_e\cap R_f \stackrel{\mu-a.e.}{=}\emptyset$ if $e \neq f \in \mathcal{G}^1$;
\item $D_\emptyset=\emptyset; D_A \cap D_B\stackrel{\mu-a.e.}{=} D_{A \cap B}; D_A \cup D_B\stackrel{\mu-a.e.}{=} D_{A \cup B}$ for all $A, B \in \mathcal{G}^0$;
\item $R_e\stackrel{\mu-a.e.}{\subseteq}D_{s(e)}$ for all $e\in \mathcal{G}^1$;
\item\label{D_v=cup_{e in s^{-1}(v)}R_e} $D_v\stackrel{\mu-a.e.}{=} \bigcup_{e \in s^{-1}(v)}R_e$ if $0 <\vert s^{-1}(v) \vert<\infty$; and
\item for each $e\in \mathcal{G}^1$, there exist two measurable maps $f_e:D_{r(e)}\rightarrow R_e$ and $f_e^{-1}:R_e \rightarrow D_{r(e)}$ such that $f_e\circ f_e^{-1}\stackrel{\mu-a.e.}{=}\id_{R_e}, f_e^{-1}\circ f_e\stackrel{\mu-a.e.}{=}\id_{D_{r(e)}}$, the pushforward measure $\mu \circ f_e$, of $f_e^{-1}$ in $D_{r(e)}$, is absolutely continuous with respect to $\mu$ in $D_{r(e)}$, and the pushforward measure $\mu \circ f_e^{-1}$, of $f_e$ in $R_e$, is absolutely continuous with respect to $\mu$ in $R_e$. Denote the Radon-Nikodym derivative $d(\mu \circ f_e)/d\mu$ by $\Phi_{f_e}$ and  the Radon-Nikodym derivative $d(\mu\circ f_e^{-1} )/d\mu$ by $\Phi_{f_e^{-1}}$.
\end{enumerate}
We call $\{R_e,D_A,f_e\}_{e \in \mathcal{G}^1,A \in \mathcal{G}^0}$ a \emph{$\mathcal{G}$-branching system} on $(X,\mu)$.
\end{definicao}

It follows that $\mu$ in $D_{r(e)}$ is absolutely continuous with respect to $\mu \circ f_e, \mu$ in $R_e$ is absolutely continuous with respect to $\mu \circ f_e^{-1}, \Phi_{f_e}>0$ $\mu$-a.e. in $D_{r(e)}, \Phi_{f_e^{-1}}>0$ $\mu$-a.e. in $R_e$, and $\Phi_{f_e}(x)\Phi_{f_e^{-1}}(f_e(x))=1$ $\mu-a.e.$ in $D_{r(e)}$. For $\alpha \in \mathcal{G}^*\setminus \mathcal{G}^0$, define $f_\alpha:=f_{\alpha_1} \circ\dots\circ f_{\alpha_n}$, and define $f_\alpha^{-1}:=f_{\alpha_n}^{-1} \circ\dots\circ f_{\alpha_1}^{-1}$. It is straightforward to see that $\mu \circ f_{\alpha}$ in $D_{r(\alpha_n)}$ is absolutely continuous with respect to $\mu$ in $D_{r(\alpha_n)}$, and $\mu \circ f_{\alpha}^{-1}$ in $R_{\alpha_1}$ is absolutely continuous with respect to $\mu$ in $R_{\alpha_1}$. Denote the Radon-Nikodym derivative $d(\mu \circ f_{\alpha})/d\mu$ by $\Phi_{f_\alpha}$, and denote the Radon-Nikodym derivative $d(\mu \circ f_{\alpha}^{-1})/d\mu$ by $\Phi_{f_\alpha^{-1}}$.

\begin{teorema}\label{existenceofabranchingsystem}
Let $\mathcal{G}$ be an ultragraph. Then there exists a $\mathcal{G}$-branching system.
\end{teorema}

\demo
Let $X:=\mathbb{R}$ and let $\mu$ be the Lebesgue measure on all Borel sets of $\mathbb{R}$. We enumerate the set $\mathcal{G}^1=\{e_i\}_{i \geq 1}$ and the set of sinks $G_{\mathrm{sink}}^0=\{v_i:s^{-1}(v_i)=\emptyset\}_{i \geq 1}$. For each $i \geq 1$, define $R_{e_i}:=[i-1,i]$, and define $D_{v_i}:=[-i,1-i]$. For $v \in G^0$ with $s^{-1}(v) \neq \emptyset$, define $D_v:=\cup_{e \in s^{-1}(v)}R_e$. Define $D_\emptyset=\emptyset$. For $A \neq \emptyset \in \mathcal{G}^0$, define $D_A:=\cup_{v \in A}D_v$.

It is easy to see that $\{R_e,D_A\}_{e \in \mathcal{G}^1,A \in \mathcal{G}^0}$ satisfies Condition~(\ref{R_e cap R_f =emptyset if e neq f})--(\ref{D_v=cup_{e in s^{-1}(v)}R_e}) of Definition~\ref{branchsystem}.

Fix $e \in \mathcal{G}^1$. We prove the existence of $f_e, f_e^{-1}, \Phi_{f_e}, \Phi_{f_e^{-1}}$. Write $R_e=[n,n+1]$ for some $n \geq 0$. Suppose that $D_{r(e)}=\bigcup_{i=1}^{m}[n_i,n_i+1]$, where $n_i <n_{i+1}$. For each $i$, let $F_i:[n_i,n_i+1] \to [n+(i-1)/m,n+i/m]$ be an arbitrary increasing bijection in $C^1([n_i,n_i+1])$. Piecing together $F_i$'s yields $f_e$, and piecing together $F_i^{-1}$'s yields $f_e^{-1}$. The existence of $\Phi_{f_e}, \Phi_{f_e^{-1}}$ follows easily.

Suppose that $D_{r(e)}=\bigcup_{i=1}^{\infty}[n_i,n_i+1]$, where $n_i <n_{i+1}$. For each $i$, let $F_i:[n_i,n_i+1] \to [n+1-(1/2)^{i-1},n+1-(1/2)^i]$ be an arbitrary increasing bijection in $C^1([n_i,n_i+1])$. Piecing together $F_i$'s yields $f_e$. Piecing together $F_i^{-1}$'s and giving an arbitrary value at $n+1$ yield $f_e^{-1}$. The existence of $\Phi_{f_e}, \Phi_{f_e^{-1}}$ follows easily. So we are done.
\fim

In \cite{Tomforde:IUMJ03} it is shown that there exists ultragraph C*-algebras that are neither Exel-Laca nor graph C*-algebras. In particular, the following example is considered:

\begin{exemplo} Let $\mathcal{G}$ be the ultragraph where $\mathcal{G}^1=\{e_i,g_i\}_{i\in \N}$, $G^0=\{w\}\cup\{v_i\}_{i\in \N}$ and with the following range and source maps: $r(g_i)=G^0\setminus\{w\}$ for each $i\in \N$, $r(e_i)=\{v_i,v_4,v_5,v_6,...\}$ for each $1\leq i\leq 3$, $r(e_i)=\{v_i,v_{i-3}\}$ for each $i\geq 4$, $s(g_i)=w$ for each $i\in \N$ and $s(e_i)=v_i$ for each $i\in \N$.
\end{exemplo}

 Since the ultragraph C*-algebra associated to the ultragraph above is not an Exel-Laca nor a graph C*-algebra, it is interesting to construct a branching system associated to $\mathcal{G}$. We will define a concrete $\mathcal{G}-$branching system in $\R$, with Lebesgue measure.

Define, for each $i\in \N$, $R_{e_i}=[i-1,i)$, $R_{g_i}=[-i,-i+1)$, $D_{v_i}=[i-1,i)$ and $D_w=(-\infty,0)$. Now, defining $D_A=\bigcup\limits_{u\in A}D_u$, for each $A\in \mathcal{G}^0$, we obtain that: $D_{r(e_i)}=[i-1,i)\cup [3,\infty)$ for each $1\leq i\leq 3$, $D_{r(e_i)}=[i-4,i-3)\cup [i-1,i)$ for each $i\geq 4$, and $D_{r(g_i)}=[0,\infty)$ for each $i\in \N$. In the next figure we show graphically an example of maps $f_{e_i}^{-1}$ and $f_{g_i}^{-1}$.
\\
\\
\centerline{
\setlength{\unitlength}{1cm}
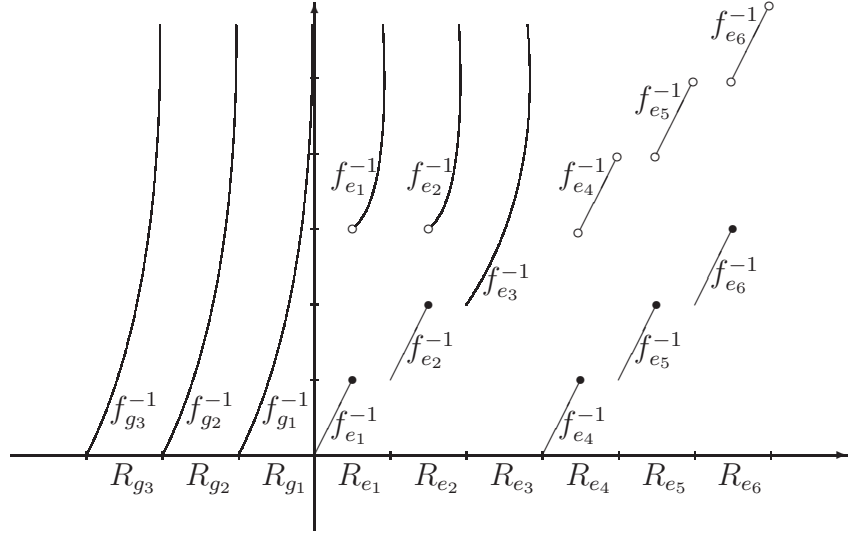
\begin{picture}(-3,0)
\put(-7,-7){\vector(1,0){11}}
%\put(3.75,-7,104){$>$}
\put(-3,-8){\vector(0,1){7}}
\put(-3.05,-6){\line(1,0){0.1}}
\put(-3.05,-5){\line(1,0){0.1}}
\put(-3.05,-4){\line(1,0){0.1}}
\put(-3.05,-3){\line(1,0){0.1}}
\put(-3.05,-2){\line(1,0){0.1}}
\put(-6,-7,05){\line(0,1){0.1}}
\put(-5,-7,05){\line(0,1){0.1}}
\put(-4,-7,05){\line(0,1){0.1}}
\put(-2,-7,05){\line(0,1){0.1}}
\put(-1,-7,05){\line(0,1){0.1}}
\put(0,-7,05){\line(0,1){0.1}}
\put(1,-7,05){\line(0,1){0.1}}
\put(2,-7,05){\line(0,1){0.1}}
\put(3,-7,05){\line(0,1){0.1}}
\put(-3.7,-7.4){$R_{g_1}$}
\put(-4.7,-7.4){$R_{g_2}$}
\put(-5.7,-7.4){$R_{g_3}$}
\put(-2.7,-7.4){$R_{e_1}$}
\put(-1.7,-7.4){$R_{e_2}$}
\put(-0.7,-7.4){$R_{e_3}$}
\put(0.3,-7.4){$R_{e_4}$}
\put(1.3,-7.4){$R_{e_5}$}
\put(2.3,-7.4){$R_{e_6}$}
\qbezier(-6,-7)(-5,-5)(-5.03,-1.3)
\qbezier(-5,-7)(-4,-5)(-4.03,-1.3)
\qbezier(-4,-7)(-3,-5)(-3.03,-1.3)
\put(-5.7,-6.5){$f_{g_3}^{-1}$}
\put(-4.7,-6.5){$f_{g_2}^{-1}$}
\put(-3.7,-6.5){$f_{g_1}^{-1}$}
\put(-3,-7){\line(1,2){0.5}}
\put(-2.5,-6){\circle*{0.1}}
\qbezier(-2.46,-3.96)(-2,-3.5)(-2.1,-1.3)
\put(-2.5,-4){\circle{0.1}}
\put(-2,-6){\line(1,2){0.5}}
\put(-1.5,-5){\circle*{0.1}}
\qbezier(-1.46,-3.96)(-1,-3.5)(-1.1,-1.3)
\put(-1.5,-4){\circle{0.1}}
\qbezier(-1,-5)(-0,-3.5)(-0.2,-1.3)
\put(0,-7){\line(1,2){0.5}}
\put(0.5,-6){\circle*{0.1}}
\put(0.5,-4){\line(1,2){0.45}}
\put(0.98,-3.05){\circle{0.1}}
\put(0.47,-4,06){\circle{0.1}}
\put(1,-6){\line(1,2){0.5}}
\put(1.5,-5){\circle*{0.1}}
\put(1.5,-3){\line(1,2){0.45}}
\put(1.48,-3.05){\circle{0.1}}
\put(1.98,-2.05){\circle{0.1}}
\put(2,-5){\line(1,2){0.5}}
\put(2.5,-4){\circle*{0.1}}
\put(2.5,-2){\line(1,2){0.45}}
\put(2.48,-2.05){\circle{0.1}}
\put(2.98,-1.05){\circle{0.1}}
\put(-2.8,-6.7){$f_{e_1}^{-1}$}
\put(-2.8,-3.4){$f_{e_1}^{-1}$}
\put(-1.8,-5.7){$f_{e_2}^{-1}$}
\put(-1.8,-3.4){$f_{e_2}^{-1}$}
\put(-0.8,-4.8){$f_{e_3}^{-1}$}
\put(0.2,-6.7){$f_{e_4}^{-1}$}
\put(0.2,-3.4){$f_{e_4}^{-1}$}
\put(1.2,-5.7){$f_{e_5}^{-1}$}
\put(1.2,-2.4){$f_{e_5}^{-1}$}
\put(2.2,-4.7){$f_{e_6}^{-1}$}
\put(2.2,-1.4){$f_{e_6}^{-1}$}
\end{picture}}
\vspace{8cm}

\begin{rmk} Notice that in the branching system above we have enumerated the edges of $\mathcal{G}$ differently from what we did in Theorem \ref{existenceofabranchingsystem}, but we kept the main idea of how to define the measurable sets $R_{e}$ and $D_A$.
\end{rmk}

\begin{rmk} We will see, Theorem \ref{repinducedbybranchingsystems}, that this branching system induces a representation of $C^*(\mathcal{G})$ in $B(\mathcal{L}^2(\R))$ and, since $\mathcal{G}$ satisfies condition~(L), this representation is faithful (see Theorem \cite[Theorem~6.7]{Tomforde:JOT03}).
\end{rmk}

%=======================================================================================
\section{Representations of Ultragraph C*-algebras on $\mathcal{L}^2(X,\mu)$ via Branching Systems}
%=======================================================================================

Next we show how to obtain a representation of an ultragraph C*-algebra from a given Branching System.

Let $\mathcal{G}$ be an ultragraph and let $\{R_e,D_A,f_e\}_{e\in \mathcal{G}, A\in \mathcal{G}^0}$ be a $\mathcal{G}$-branching system on $(X,\mu)$.
Since the domain of $\Phi_{f_e^{-1}}$ and the domain of $ f_e^{-1}$ are $R_e$, we can also regard them as measurable maps on $X$ by simply extending then with value zero out of $R_e$, and so, for each $\phi\in \mathcal{L}^2(X,\mu)$, we can consider the function $\Phi_{f_e^{-1}}^{1/2}\cdot( \phi \circ f_e^{-1})$. Also by extending $f_e$ and $\Phi_{f_e}$ by zero out of $D_{r(e)}$ we get the function $\Phi_{f_e}^{1/2}\cdot( \phi \circ f_e)$.

\begin{teorema}\label{repinducedbybranchingsystems}
Let $\mathcal{G}$ be an ultragraph and let $\{R_e,D_A,f_e\}_{e\in \mathcal{G}^1, A\in \mathcal{G}^0}$ be a $\mathcal{G}$-branching system on a measure space $(X,\mu)$. Then there exists an unique representation $\pi:C^*(\mathcal{G}) \to B(\mathcal{L}^2(X,\mu))$ such that $\pi(s_e)(\phi)=\Phi_{f_e^{-1}}^{1/2}\cdot( \phi \circ f_e^{-1})$ and $\pi(p_A)(\phi)=\chi_{D_A}\phi$, for all $e \in \mathcal{G}^1, A \in \mathcal{G}^0$ and $\phi \in \mathcal{L}^2(X,\mu)$.
\end{teorema}

\demo
It is straightforward to check that $\{\pi(p_A)\}_{A \in \mathcal{G}^0}$ is a family of projections in $B(\mathcal{L}^2(X,\mu))$ that satisfies Condition~(\ref{p_Ap_B=p_{A cap B}}) of Definition~\ref{defofultragraphalgebra}.

Fix $e \in \mathcal{G}^1$. For $\phi \in \mathcal{L}^2(X,\mu)$, we have
\[
\int \vert \Phi_{f_e^{-1}}^{1/2}\cdot(\phi \circ f_e^{-1})\vert^2\, \mathrm{d}\mu=\int \vert \phi\circ f_e^{-1}\vert^2\, \mathrm{d}(\mu\circ f_e^{-1})=\int_{D_{r(e)}} \vert \phi \vert^2\, \mathrm{d}\mu<\infty.
\]
So $\pi(s_e)(\phi) \in \mathcal{L}^2(X,\mu)$. Define
\begin{equation}
\pi(s_e)^*(\phi):= \Phi_{f_e}^{1/2}\cdot( \phi \circ f_e).
\end{equation}
Similarly $\pi(s_e)^*$ is an operator on $\mathcal{L}^2(X,\mu)$. For $\phi,\eta \in \mathcal{L}^2(X,\mu)$, we have
\begin{align*}
\langle \pi(s_e)^*(\phi),\eta \rangle&=\int \Phi_{f_e}^{1/2}\cdot(\phi \circ f_e)\cdot\overline{\eta}\, \mathrm{d}\mu
\\&=\int \Phi_{f_e}^{-1/2}\cdot(\phi \circ f_e)\cdot\overline{\eta}\, \mathrm{d}(\mu\circ f_e)
\\&=\int (\Phi_{f_e}^{-1/2}\circ f_e^{-1})\cdot\phi\cdot(\overline{\eta \circ f_e^{-1}})\, \mathrm{d}\mu
\\&=\int \phi \cdot \Phi_{f_e^{-1}}^{1/2}\cdot\overline{\eta \circ f_e^{-1}}\, \mathrm{d}\mu
\\&=\langle \phi,\pi(s_e)(\eta)\rangle.
\end{align*}
So $\pi(s_e)^*$ is the adjoint of $\pi(s_e)$.

For $e \in \mathcal{G}^1, \phi \in \mathcal{L}^2(X,\mu)$, we have that
\begin{equation}
\pi(s_e)\pi(s_e)^*(\phi)\stackrel{\mu-a.e.}{=} \chi_{R_e}\phi.
\end{equation}
Since $R_e \cap R_f=\emptyset$ if $e \neq f \in \mathcal{G}^1$, we get a family of partial isometries with orthogonal ranges $\{\pi(s_e)\}_{e \in \mathcal{G}^1}$. We also get that $\pi(s_e)\pi(s_e)^*\leq \pi(p_{s(e)})$ because $R_e \stackrel{\mu-a.e.}{\subset} D_{s(e)}$.

Condition~(\ref{s_e^*s_e=p_{r(e)}}) of Definition~\ref{defofultragraphalgebra} follows easily.

Finally we check Condition~(\ref{CK-condition}) of Definition~\ref{defofultragraphalgebra}. Fix $v \in G^0$ with $0<\vert s^{-1}(v)\vert<\infty$. For $\phi \in \mathcal{L}^2(X,\mu)$, since $\bigcup_{e \in s^{-1}(v)}R_e\stackrel{\mu-a.e.}{=}D_{s(e)}$, we have that
\begin{align*}
\sum_{e\in s^{-1}(v)}\pi(s_e)\pi(s_e)^*(\phi)=\sum_{e \in s^{-1}(v)}\chi_{R_e}\phi\stackrel{\mu-a.e.}{=}\pi(p_{v})(\phi)
\end{align*}
and therefore we are done. \fim

\begin{rmk}
For any $\phi \in L^2(X,\mu)$, we have $\pi(s_\alpha)(\phi)=\Phi_{f_{\alpha}^{-1}}^{1/2}\phi \circ f_{\alpha}^{-1}$, and $\pi(s_\alpha)^*(\phi)=\Phi_{f_{\alpha}}^{1/2}\phi \circ f_{\alpha}$.
\end{rmk}

\begin{corolario}\label{existence of a representation on L^2(R) of C^*(G)}
Let $\mathcal{G}$ be an ultragraph. Then there exists a $\mathcal{G}$-branching system $\{R_e,D_A,f_e\}_{e\in \mathcal{G}^1,A\in \mathcal{G}^0}$ on $(\mathbb{R},\mu)$, where $\mu$ is the Lebesgue measure on all Borel sets of $\mathbb{R}$, and there exists a unique representation $\pi:C^*(\mathcal{G}) \to B(\mathcal{L}^2(\mathbb{R},\mu))$, such that $\pi(s_e)(\phi)=\Phi_{f_e^{-1}}^{1/2}( \phi \circ f_e^{-1})$ and $\pi(p_A)(\phi)=\chi_{A}\phi$, for all $e \in \mathcal{G}^1, A \in \mathcal{G}^0$ and $\phi \in \mathcal{L}^2(\mathbb{R},\mu)$.
\end{corolario}

\demo
It follows immediately from Theorem~\ref{existenceofabranchingsystem} and Theorem~\ref{repinducedbybranchingsystems}.
\fim

%========================================================================================================================
\section{Nonsingular Branching Systems of Ultragraphs}\label{PFoperator}
%========================================================================================================================

The Perron-Frobenius operator (or Frobenius-Perron operator, or transfer operator) of ergodic theory is used, among other things, to study invariant measures of non invertible transformations (see \cite{LY} for example).

In this section we describe the Perron-Frobenius operator in terms of the representations of ultragraph C*-algebras introduced in the previous section.

\begin{definicao}\label{defofnonsingope}
Let $F$ be a measurable map on a measure space $(X,\mu)$. Then $F$ is called \emph{nonsingular} if for any measurable set $S \subset X, \mu(S)=0$ implies that $\mu(F^{-1}(S))=0$. The unique operator $P_F \in B(\mathcal{L}^1(X,\mu))$ such that $\int_S P_F \phi\, \mathrm{d}\mu=\int_{F^{-1}(S)}\phi \, \mathrm{d}\mu$, for all $\phi \in \mathcal{L}^1(X,\mu)$, for all measurable set $S \subset X$,  is called the \emph{Perron-Frobenius operator} corresponding to $F$.
\end{definicao}

\begin{definicao}\label{defofnonsingbranchingsys}
Let $\mathcal{G}$ be an ultragraph and let $\{R_e,D_A,f_e\}_{e\in \mathcal{G}^1,A\in\mathcal{G}^0}$ be a $\mathcal{G}$-branching system on a measure space $(X,\mu)$. The $\mathcal{G}$-branching system is called \emph{nonsingular} if there exists a nonsingular map $F:X\to X$ such that $F \vert_{R_e}\stackrel{\mu-a.e.}{=}f_e^{-1}$ for all $e \in \mathcal{G}^1$.
\end{definicao}

\begin{teorema}\label{branching sys is nonsingular}
Let $\mathcal{G}$ be an ultragraph. Then any $\mathcal{G}$-branching system $\{R_e,$ $D_A, f_e\}_{e\in \mathcal{G}^1,A\in\mathcal{G}^0}$ on a measure space $(X,\mu)$ is nonsingular.
\end{teorema}

\demo
Fix $e \in \mathcal{G}^1$. Define a measurable map $F_e: R_e \setminus (\bigcup_{e' \neq e}R_{e'}) \to X$ by $F_e(x):=f_e^{-1}(x)$. Take an arbitrary constant map $g:\bigcup_{e \neq e'}(R_e \cap R_{e'}) \to X$. Let $h:X \setminus \bigcup_{e \in \mathcal{G}^1}R_e \to X$ be the inclusion map. Then we get a measurable map $F:X \to X$ by piecing together $F_e$'s, $g$ and $h$. It is obvious that $F \vert_{R_e}\stackrel{\mu-a.e.}{=}f_e^{-1}$ for all $e \in \mathcal{G}^1$.

Now we check that $F$ is nonsingular. Fix a measurable set $S$ such that $\mu(S)=0$. Then $F^{-1}(S)$ is the disjoint union of all the $F_e^{-1}(S), g^{-1}(S)$ and $h^{-1}(S)$. We observe that $\mu(g^{-1}(S))=0$, since $\mu(\bigcup_{e \neq e'}(R_e \cap R_{e'}))=0$, and that $\mu(h^{-1}(S))=0$ because $h^{-1}(S) \subset S$. For $e \in \mathcal{G}^1$, we have $F_e^{-1}(S)\stackrel{\mu-a.e.}{=}f_e(S \cap D_{r(e)})$. Since $\mu \circ f_e$ is absolutely continuous with respect to $\mu$, we get $\mu(F_e^{-1}(S))=0$. So $\mu(F^{-1}(S))=0$ and we are done.
\fim

\begin{corolario}\label{existence of a nonsing mathcal{G}-branching system}
Let $\mathcal{G}$ be an ultragraph. Then there exists a nonsingular $\mathcal{G}$-branching system.
\end{corolario}

\demo
It follows immediately from Theorem~\ref{existenceofabranchingsystem} and Theorem~\ref{branching sys is nonsingular}.
\fim

Next we show the relation between the Perron-Frobenius operator and representations of ultragraph C*-algebras.

\begin{teorema}\label{Perron-Frob ope}
Let $\mathcal{G}$ be an ultragraph graph, $\{R_e,D_A,f_e\}_{e\in \mathcal{G}^1,A\in\mathcal{G}^0}$ be a $\mathcal{G}$-branching system on a measure space $(X,\mu)$. Suppose that $F$ is a nonsingular map on $X$ such that $F \vert_{R_e}\stackrel{\mu-a.e.}{=}f_e^{-1}$, for all $e \in \mathcal{G}^1$, and let $\pi:C^*(\mathcal{G}) \to B(\mathcal{L}^2(X,\mu))$ be the representation from Theorem~\ref{repinducedbybranchingsystems}. Then for $\phi\in\mathcal{L}^2(X,\mu), \eta \in \mathcal{L}^1(X,\mu)$, we have that
\begin{enumerate}
\item if $supp(\phi)\stackrel{\mu-a.e.}\subset \bigcup_{i=1}^{n} R_{e_i}$, then $P_F(\phi^2)\stackrel{\mu-a.e.}=\sum_{i=1}^{n}(\pi(s_{e_i}^*)\phi)^2$;
\item if $supp(\phi)\stackrel{\mu-a.e.}\subset \bigcup_{i=1}^\infty R_{e_i}$, then $P_F(\phi^2)\stackrel{\mu-a.e.}=\lim_{n \to\infty}\sum_{i=1}^n(\pi(s_{e_i}^*)\phi)^2$ under the $\mathcal{L}^1(X,\mu)$-norm; and
\item if we write $\eta=\eta_1-\eta_2+i(\eta_3-\eta_4)$, where each $\eta_i$ is nonnegative $\mu$-a.e., and $supp(\eta_i)\stackrel{\mu-a.e.}\subset \bigcup_{i=1}^\infty R_{e_i}$, then
\begin{align*}
P_F(\eta)=\lim_{n \to\infty}\sum_{i=1}^n\Big((\pi(s_{e_i}^*)\sqrt{\eta_1})^2-(\pi(s_{e_i}^*)\sqrt{\eta_2})^2+&i(\pi(s_{e_i}^*)\sqrt{\eta_3})^2-\\&i(\pi(s_{e_i}^*)\sqrt{\eta_4})^2\Big)
\end{align*}
under the $\mathcal{L}^1(X,\mu)$-norm.
\end{enumerate}
\end{teorema}

\demo
We prove the first statement. Fix a measurable set $S \subset X$. Then
\begin{align*}
\int_S\sum_{i=1}^{n}(\pi(s_{e_i}^*)\phi)^2\, \mathrm{d}\mu&=\sum_{i=1}^{n}\int_S \Phi_{f_{e_i}}\cdot(\phi^2\circ f_{e_i})\, \mathrm{d}\mu
\\&=\sum_{i=1}^{n}\int_S \phi^2\circ f_{e_i}\, \mathrm{d}(\mu\circ f_{e_i})
\\&=\sum_{i=1}^{n}\int_{f_{e_i}(S \cap D_{r(e_i)})} \phi^2\, \mathrm{d}\mu
\\&=\sum_{i=1}^{n}\int_{F^{-1}(S) \cap R_{e_i}} \phi^2\, \mathrm{d}\mu
\\&=\int_{F^{-1}(S)} \phi^2\, \mathrm{d}\mu
\\&=\int_{S} P_F(\phi^2)\, \mathrm{d}\mu.
\end{align*}
So $P_F(\phi^2)=\sum_{i=1}^{n}(\pi(s_{e_i}^*)\phi)^2$.

Now we prove the second statement. By the dominated convergence theorem, we have that $(\phi \chi_{(\bigcup_{i=1}^{n}R_{e_i})})^2 \to\phi^2$ under the $\mathcal{L}^1(X,\mu)$-norm. By continuity of $P_F$, and by the first statement, we have that
\begin{align*}
P_F(\phi^2)&=\lim_{n \to \infty}P_F((\phi \chi_{(\bigcup_{i=1}^{n}R_{e_i})})^2)
\\&=\lim_{n \to \infty}\sum_{i=1}^{n}(\pi(s_{e_i}^*)(\phi \chi_{(\bigcup_{j=1}^{n}R_{e_j})}))^2
\\&=\lim_{n \to \infty}\sum_{i=1}^{n}(\pi(s_{e_i}^*)\phi)^2.
\end{align*}

Straightforward calculations yield the third statement.
\fim

%================================================================================================================================
\section{Permutative Representations of Ultragraph C*-algebras}
%================================================================================================================================

Permutative representations of the Cuntz algebra, and their relations to iterated function systems, were studied in the seminal work of Bratteli and Jorgensen, see \cite{BJ}. Many authors have generalized and further studied these representations. For example, motivated by the work of Kawamura in \cite{Kawamura}, Lawson studied primitive partial permutation representations of polycyclic monoids and their relations to branching function systems (see \cite{Lawson}).

In this section we define permutative representations of ultragraph C*-algebras and show that every permutative representation is unitary equivalent to a representation arising from a branching system. We then proceed to show that for a class of ultragraphs every representation is permutative.

Let $\mathcal{G}$ be an ultragraph and let $\pi:C^*(\mathcal{G}) \to B(H)$ be a representation. For $e \in \mathcal{G}^1$, define a closed subspace of $H$ by $H_e:=\pi(s_es_e^*)(H)$. For $A \in \mathcal{G}^0$, define a closed subspace of $H$ by $H_A:=\pi(p_A)(H)$. Then we have the following.
\begin{enumerate}
\item $H_A\cap H_B=H_{A\cap B}$ for $A,B \in \mathcal{G}^0$;
\item $H_e \perp H_f$ if $e \neq f$;
\item $\pi(s_e) \vert_{H_{r(e)}}:H_{r(e)} \to H_e$ is an isomorphism of Hilbert spaces;
\item for $v \in G^0$ with $s^{-1}(v) \neq \emptyset$, we have $H_v=\Big(\bigoplus\limits_{e \in s^{-1}(v)}H_e\Big) \oplus V_v$;
\item for $v \in G^0$ with $0 <\vert s^{-1}(v) \vert <\infty$, we have $V_v=0$; and
\item $H=\Big(\bigoplus\limits_{v \in G^0}H_v\Big) \oplus V$.
\end{enumerate}

In \cite[Page~119]{MR2848777} permutative representations of graph algebras were defined (though the authors did not explicitly use the term permutative representation). Below we generalize the definition in \cite{MR2848777} to ultragraphs.

\begin{definicao}\label{define permutative rep}
Let $\mathcal{G}$ be an ultragraph and let $\pi:C^*(\mathcal{G}) \to B(H)$ be a representation. Then $\pi$ is said to be \emph{permutative} if there exist orthonormal bases $B$ of $H; B_v$ of $H_v$ for all $v \in G^0; B_{r(e)}$ of $H_{r(e)}$ for all $e \in \mathcal{G}^1$; and $B_e$ of $H_e$ for all $e \in \mathcal{G}^1$ such that
\begin{enumerate}
\item $B_{r(e)} \subset B$ for all $e \in \mathcal{G}^1$ and $B_v\subseteq B$ for all $v\in G^0$;
\item $B_v \subset B_{r(e)}$ for all $v \in G^0, e \in \mathcal{G}^1$ with $v \in r(e)$;
\item $B_v \supset \bigcup\limits_{e \in s^{-1}(v)}B_e$ for all $v \in G^0$; and
\item $\pi(s_e)(B_{r(e)})=B_e$ for all $e \in \mathcal{G}^1$ (\emph{B2B}).
\end{enumerate}
\end{definicao}

\begin{rmk}
The $B2B$ condition is equivalent to $\pi(s_e^*)(B_e)=B_{r(e)}$. for all $e \in \mathcal{G}^1$.
\end{rmk}

\begin{lema}\label{Define B_A from perm rep}
Let $\mathcal{G}$ be an ultragraph, $\pi:C^*(\mathcal{G}) \to B(H)$ be a permutative representation and let $B$ be an orthonormal basis of $H$ satisfying the conditions of Definition~\ref{define permutative rep}. For $A\in \mathcal{G}^0$, following Lemma~\ref{concrete des of G^0}, let $A:=\left(\bigcap\limits_{e\in X_1}r(e)\right)\cup...\cup \left(\bigcap\limits_{e\in X_n}r(e)\right)\cup F$ and define
\[
B_A:=\left(\bigcap\limits_{e\in X_1}B_{r(e)}\right)\cup...\cup \left(\bigcap\limits_{e\in X_n}B_{r(e)}\right)\cup \bigcup\limits_{v\in F}B_v.
\]
Then the set $B_A$ is well-defined. Furthermore, $B_A$ is an orthonormal basis of $H_A$ and, for $A,C\in \mathcal{G}^0$, we have that $B_A\cap B_C=B_{A\cap C}$ and $B_A\cup B_C=B_{A\cup C}$.
\end{lema}

\demo
Fix $A\in \mathcal{G}^0$ and write $A=\left(\bigcap\limits_{e\in X_1}r(e)\right)\cup...\cup \left(\bigcap\limits_{e\in X_n}r(e)\right)\cup F$. Note that $B_A\subseteq B$ since $B_{r(e)}\subseteq B$ and $B_v\subseteq B$ for all $e \in \mathcal{G}^1$ and all $v \in G^0$. First we prove that $B_A$ is an orthonormal basis of $H_A:=\pi(p_A)(H)$. Fix $h\in B_A$. If $h\in B_v$ for some $v\in F$ then $\pi(p_A)(\pi(p_v)(h))=\pi(p_{A\cap v})(h)=\pi(p_v)(h)=h \in H_A$. If $h\in \bigcap\limits_{e\in X_i}B_{r(e)}\subseteq \bigcap\limits_{e\in X_i}H_{r(e)}$ for some $i$, then since $\bigcap\limits_{e\in X_i} r(e)\subseteq A$ we get $\pi(p_A)\pi(p_{\bigcap\limits_{e\in X_i}r(e)})(h)=\pi(p_{(A\bigcap\limits_{e\in X_i}r(e))})(h)=\pi(p_{\bigcap\limits_{e\in X_i}r(e)})(h)=h\in H_A$. Therefore, $B_A\subseteq H_A$. Let $V$ be the closed subspace of $H_A$ generated by $B_A$. Let $W$ be the orthogonal complement of $V$ in $H_A$, so that $H_A=V\oplus W$. Now we show that $W=0$. For each $v\in F$ we get $\pi(p_v)(W)=0$ since $\pi(p_v)(H)\subseteq V$ and $W$ is orthogonal to $V$. Similarly, since $B_{r(e)} \subset B$ for all $e$, we have that $\pi(p_{\bigcap\limits_{e\in X_i}r(e)})(W)=0$ for each $i$. We deduce that $\pi(p_A)(W)=0$. So $W=0$ because $W=\pi(p_A)(W)$. It follows that the set $B_A$ is an orthonormal basis of $H_A$ and $B_A\subseteq B$. Since there is only one subset of $B$ which is an orthonormal basis of $H_A$, it follows that $B_A$ is well-defined (notice that since $A \in \mathcal{G}^0$ can be described in more then one way using Lemma~\ref{concrete des of G^0}, we had to show that $B_A$ is well defined).

Now we show that for each $A,C\in \mathcal{G}^0, B_{A\cup C}=B_A\cup B_C$ and $B_{A\cap C}=B_A\cap B_C$.

It is straightforward to see that $B_A \cup B_C \subset B_{A \cup C}$ because $\pi(p_A)H$, $\pi(p_C)H \subset \pi(p_{A \cup C})H$. Conversely, for $h \in B_{A \cup C}$, suppose that $h \notin B_A \cup B_C$. Then $h=\pi(p_{A \cup C})h=\pi(p_A)h+\pi(p_C)h-\pi(p_{A\cap C})h=0$ which is a contradiction. So $B_{A\cup C}=B_A\cup B_C$.

As before, it is straightforward to see that $B_{A \cap C} \subset B_A \cap B_C$ because $\pi(p_{A \cap C})H \subset \pi(p_A)H, \pi(p_C)H$. Conversely, for $h \in B_A \cap B_C$, we have that $h=\pi(p_A)\pi(p_C)h=\pi(p_{A \cap C})h$. So $h \in B_{A \cap C}$ and hence $B_{A \cap C}=B_A \cap B_C$. \fim

The following theorem is a generalization of \cite[Theorem~2.1]{MR2848777}.

\begin{teorema}\label{unitequiv}
Let $\mathcal{G}$ be an ultragraph, $\pi:C^*(\mathcal{G}) \to B(H)$ be a permutative representation and let $B$ be an orthonormal basis of $H$ satisfying conditions of Definition~\ref{define permutative rep}. Then $\pi$ is unitarily equivalent to a representation of $C^*(\mathcal{G})$ on $l^2(\mathbb{N})$ which is induced from a $\mathcal{G}$-branching system on $\mathbb{N}$ (the measure on $\mathbb{N}$ is the counting measure).
\end{teorema}

\demo
Write $B:=\{h_n\}_{n \in \mathbb{N}}$. Then we get a unitary $U:H \to l^2(\mathbb{N})$ such that $U(h_n)=\delta_n$, where $\delta_n$ is the point mass function. For $e \in \mathcal{G}^1$, define $R_e:=\{n \in \mathbb{N}:h_n \in B_e\}$. For $A\in \mathcal{G}^0$, define $D_A:=\{n\in \N:h_n\in B_A\}$ (see Lemma~\ref{Define B_A from perm rep}). Then for $e \in \mathcal{G}^1$, the (B2B) condition yields a bijection $f_e:D_{r(e)} \to R_e$ such that $h_{f_e(n)}=\pi(s_e)(h_n)$ and $h_{f_e^{-1}(n)}=\pi(s_e^*)(h_n)$. It is straightforward to check that $\{R_e,D_A,f_e\}_{e \in \mathcal{G}^1, A \in \mathcal{G}^0}$ is a $\mathcal{G}$-branching system on $\mathbb{N}$. By Theorem~\ref{repinducedbybranchingsystems}, there exists a unique representation $\rho:C^*(\mathcal{G}) \to B(l^2(\mathbb{N}))$ induced from the branching system $\{R_e,D_A,f_e\}$. Using techniques similar to the graph case, see \cite[Theorem~2.1]{MR2848777}, one can show that $\pi$ is unitarily equivalent to $\rho$ via $U$. \fim

We now proceed to describe ultragraphs for which a large class of representations is permutative. Notice that while for Theorem \ref{unitequiv} above the techniques used for graphs in \cite{MR2848777} are readily applied, the same is not true for the reminder of this section.

\begin{definicao}
Let $\mathcal{G}$ be an ultragraph. An \emph{extreme vertex} is an element $A\in r(\mathcal{G}^1)\cup s(\mathcal{G}^1)$ satisfying
\begin{enumerate}
\item either $A=r(e)$ for some edge $e$ and $A\cap r(\mathcal{G}^1\setminus\{e\})=\emptyset=A\cap s(\mathcal{G}^1)$; or
\item $A=s(e)$ for some edge $e$ and $A\cap s(\mathcal{G}^1\setminus \{e\})=\emptyset=A\cap r(\mathcal{G}^1)$.
\end{enumerate}
The edge $e$ associated to an extreme vertex $A$ as above is called the \emph{extreme edge} of $A$.
\end{definicao}

For an ultragraph $\mathcal{G}$ let $X_1$ be the set of extreme vertices and let $Y_1$ be the set of extreme edges.

\begin{exemplo} Let $\mathcal{G}$ be the ultragraph of Example \ref{exeUltragrafo} (For reader convenience we draw $\mathcal{G}$ again below). We have that $\Big(\bigcup_{e \in \mathcal{G}^1}r(e)\Big) \cup s(\mathcal{G}^1)=G^0$. The extreme vertices of this ultragraph are $\{v_1\}, \{v_7\}$ and $\{v_8,v_9\}$ so that $X_1=\{v_1,v_7,\{v_8,v_9\} \}$. Furthermore, the extreme edges are $e_1, e_3$ and $e_5$, so that $Y_1=\{e_1,e_3,e_5\}$.

\vspace{1 cm}
\centerline{
\setlength{\unitlength}{1cm}
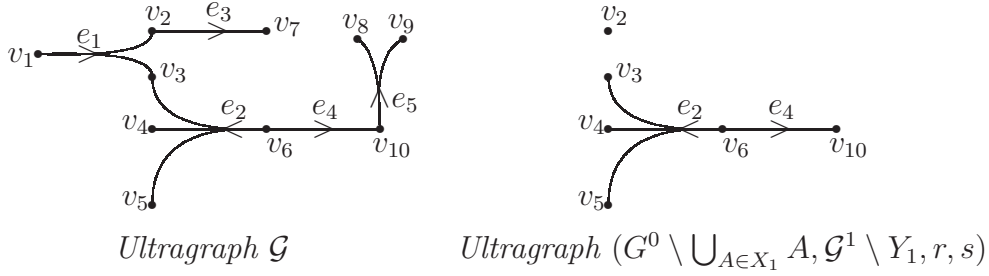
\begin{picture}(11,0)
\put(-0.5,0){\circle*{0.1}}
\put(-0.9,-0.1){$v_1$}
\qbezier(-0.5,0)(1,0.05)(1,-0.3)
\qbezier(-0.5,0)(1,-0.05)(1,0.3)
\put(1,0.3){\circle*{0.1}}
\put(0.9,0.45){$v_2$}
\put(1,-0.3){\circle*{0.1}}
\put(1.1,-0.3){$v_3$}
\put(0,-0.1){$>$}
\put(0,0.15){$e_1$}
\put(2.5,-1){\circle*{0.1}}
\put(2.5,-1.3){$v_6$}
\qbezier(2.5,-1)(1,-1.05)(1,-0.3)
\qbezier(2.5,-1)(1,-1)(1,-1)
\qbezier(2.5,-1)(1,-0.9)(1,-2)
\put(1,-1){\circle*{0.1}}
\put(0.6,-1){$v_4$}
\put(1,-2){\circle*{0.1}}
\put(0.6,-2){$v_5$}
\put(1.9,-1.1){$<$}
\put(1.9,-0.8){$e_2$}
%\put(2.5,0){\circle*{0.1}}
%\put(2.1,0){$v_7$}
\put(2.5,0.3){\circle*{0.1}}
\put(2.6,0.3){$v_7$}
\put(1,0.3){\line(1,0){1.5}}
\put(1.7,0.2){$>$}
\put(3.1,-0.8){$e_4$}
\put(4,-1){\circle*{0.1}}
\put(3.9,-1.3){$v_{10}$}
\put(2.5,-1){\line(1,0){1.5}}
\put(3.1,-1.1){$>$}
\put(1.7,0.5){$e_3$}
\qbezier(4,-1)(3.9,0)(4.3,0.2)
\qbezier(4,-1)(4,0)(3.7,0.2)
\put(3.85,-0.7){$\wedge$}
\put(4.15,-0.7){$e_5$}
\put(3.7,0.2){\circle*{0.1}}
\put(3.5,0.35){$v_8$}
\put(4.3,0.2){\circle*{0.1}}
\put(4.1,0.35){$v_9$}
\put(0.5,-2.7){Ultragraph $\mathcal{G}$}
\put(7,0.3){\circle*{0.1}}
\put(6.9,0.45){$v_2$}
\put(7,-0.3){\circle*{0.1}}
\put(7.1,-0.3){$v_3$}
\put(8.5,-1){\circle*{0.1}}
\put(8.5,-1.3){$v_6$}
\qbezier(8.5,-1)(7,-1.05)(7,-0.3)
\qbezier(8.5,-1)(7,-1)(7,-1)
\qbezier(8.5,-1)(7,-0.9)(7,-2)
\put(7,-1){\circle*{0.1}}
\put(6.6,-1){$v_4$}
\put(7,-2){\circle*{0.1}}
\put(6.6,-2){$v_5$}
\put(7.9,-1.1){$<$}
\put(7.9,-0.8){$e_2$}
\put(9.1,-0.8){$e_4$}
\put(10,-1){\circle*{0.1}}
\put(9.9,-1.3){$v_{10}$}
\put(8.5,-1){\line(1,0){1.5}}
\put(9.1,-1.1){$>$}
\put(5.0,-2.7){Ultragraph $(G^0\setminus \bigcup_{A \in X_1}A,\mathcal{G}^1\setminus Y_1,r,s)$}
\end{picture}}
\vspace{3.5cm}
Note that $v_2$ is an isolated vertex of the ultragraph $(G^0\setminus \bigcup_{A \in X_1}A, \mathcal{G}^1\setminus Y_1,r,s)$.
\end{exemplo}

Let $\mathcal{G}$ be an ultragraph. Define the set of \emph{isolated vertices} of $\mathcal{G}$ to be
\[
I_0:=\Big\{v\in G^0: v \notin \Big(\Big(\bigcup_{e \in \mathcal{G}^1}r(e)\Big) \cup s(\mathcal{G}^1)\Big)\Big\}
\]
and define the ultragraph $\mathbb{G}_0:=(G^0 \setminus I_0, \mathcal{G}^1,r,s)$. Denote by $X_1$ the set of extreme vertices of $\mathbb{G}_0$, let $\overline{X_1}=\bigcup\limits_{A\in X_1}A$, and denote by $Y_1$ the set of extreme edges of $\mathbb{G}_0$. Notice that the extreme vertices and the extreme edges of $\mathcal{G}$ and $\mathbb{G}_0$ are the same. Denote by $I_1$ the set of isolated vertices of the ultragraph $\Big(G^0 \setminus (I_0 \cup \overline{X_1}),\mathcal{G}^1\setminus Y_1,r,s\Big)$, and define
$$\mathbb{G}_1=\Big(G^0 \setminus (I_0 \cup I_1\cup \overline{X_1}),\mathcal{G}^1\setminus Y_1,r,s\Big).$$
Now, define $X_2$ and $Y_2$ as being the extreme vertices and extreme edges of the ultragraph $\mathbb{G}_1$, let $\overline{X_2}=\bigcup\limits_{A\in X_2}A$, let $I_2$ be the isolated vertices of the ultragraph $$\Big(G^0 \setminus \big(I_0 \cup I_1\cup \overline{X_1}\cup\overline{X_2}\big),\mathcal{G}^1\setminus (Y_1\cup Y_2),r,s\Big)$$ and let

$$\mathbb{G}_2=\Big(G^0 \setminus \big(I_0 \cup I_1\cup I_2\cup \overline{X_1}\cup\overline{X_2}\big),\mathcal{G}^1\setminus (Y_1\cup Y_2),r,s\Big).$$

Inductively, while $X_n\neq\emptyset$, we define the ultragraphs $\mathbb{G}_n$ and the sets $X_{n+1}$, of extreme vertices of $\mathbb{G}_n$, and $Y_{n+1}$, of extreme edges $\mathbb{G}_n$. We also define the sets $\overline{X_{n+1}}=\bigcup\limits_{A\in X_n}A$ and the set of isolated vertices $I_{n+1}$ of the ultragraph $\mathbb{G}_n$.

Notice that there is a bijective correspondence between the sets $X_n$ and $Y_n$, associating each extreme vertex $A\in X_n$ to an unique extreme edge $e\in Y_n$. For each $A\in X_n$, let $e\in Y_n$ be the (unique) edge associated to $A$. If $A=r(e)$ then $A$ is called a {\it final vertex} of $X_n$ and, if $A=s(e)$, then $A$ is called an {\it initial vertex} of $X_n$. We denote the set of initial vertices of $X_n$ by $X_n^{\mathrm{ini}}$ and the set of final vertices of $X_n$ by $X_n^{\mathrm{fin}}$.

\begin{lema}\label{l1}
Let $\mathcal{G}$ be an ultragraph. Suppose that there exists $n \geq 1$ such that $X_1,\dots,X_n \neq \emptyset$ and $\Big(\bigcup\limits_{e \in \mathcal{G}^1}r(e)\Big) \cup s(\mathcal{G}^1)=\bigcup\limits_{i=1}^{n}(\overline{X_i} \cup I_i)$. Then
\begin{enumerate}
\item for $1 \leq N\leq n, A,B \in X_N$, we have $A \cap B=\emptyset$;
\item for $1 \leq N \leq n, A \in X_N^{\mathrm{fin}}$ and for $e \in s^{-1}(A)$, we have $r(e)\in \bigcup\limits_{i=1}^{N-1}X_i^{\mathrm{fin}}$;
\item for $1 \leq N \leq n, v \in I_N$ and for $e \in s^{-1}(v)$, we have $r(e) \in \bigcup\limits_{i=1}^{N} X_i^{\mathrm{fin}}$; and
\item for $1 \leq N \leq n, v \in X_N^{\mathrm{ini}}$ and for $e \in s^{-1}(v)$, we have $r(e)\subset \bigcup_{i=1}^{n}$ $\Big(\Big(\bigcup\limits_{A \in X_i^{\mathrm{fin}}}A\Big) \cup I_i\Big) \cup \bigcup\limits_{i=N+1}^{n}X_i^{\mathrm{ini}}$.
\end{enumerate}
\end{lema}

\demo 1. It is straightforward.

2. Let $A\in X_N^{\mathrm{fin}}$ an extreme vertex of level $N$ and let $f\in Y_N$ be the extreme edge of level $N$ associated to $A$, that is, $r(f)=A$. Let $e\in s^{-1}(A)$. Note that since $A$ is an extreme vertex of level $N$ then $e\in Y_i$ for some $1\leq i\leq N-1$. Since $e\in Y_i$ then the extreme set associated to $e$ is an element of $X_i$ and, since $s(e)\in A\in X_N$, then the extreme set of $e$ is $r(e)\in X_i^{\mathrm{fin}}$.

3. Let $v\in I_N$ and let $e\in s^{-1}(v)$.
Since $v\in I_N$ then $e\in Y_i$ for some $1\leq i\leq N$. Moreover, since $s(e)$ is not an extreme edge (since $s(e)=v\in I_N$) then $r(e)$ is an extreme set. Since $e\in Y_i$ then $r(e)\in X_i^{\mathrm{fin}}$.

4. Fix $1 \leq N \leq n, v \in X_N^{\mathrm{ini}},e \in s^{-1}(v)$. Suppose that $r(e) \cap \bigcup_{i=1}^{N}X_i^{\mathrm{ini}} \neq \emptyset$, for a contradiction. Then there exist $1 \leq K\leq N, w \in r(e) \cap X_K^{\mathrm{ini}}$ such that $w \neq v$. Let $f\in Y_K$ be the associated extreme edge of $w$ with $s(f)=w$. Notice that $e \neq f$. In the ultragraph $\mathbb{G}_{K-1},v,w \in \mathbb{G}_{K-1}^0, e,f \in \mathcal{G}_{K-1}^1$. However, we get $w \notin X_K^{\mathrm{ini}}$, which is a contradiction. So $r(e) \cap \bigcup_{i=1}^{N}X_i^{\mathrm{ini}} = \emptyset$. Since $\Big(\bigcup\limits_{e \in \mathcal{G}^1}r(e)\Big) \cup s(\mathcal{G}^1)=\bigcup\limits_{i=1}^{n}(\overline{X_i} \cup I_i)$, we conclude that $r(e)\subset \bigcup_{i=1}^{n}\Big(\Big(\bigcup_{A \in X_i^{\mathrm{fin}}}A\Big) \cup I_i\Big) \cup \bigcup_{i=N+1}^{n}X_i^{\mathrm{ini}}$. \fim

The following theorem is an ultragraph version of \cite[Theorem~4.1]{MR2848777}.

\begin{teorema}
Let $\mathcal{G}$ be an ultragraph and let $\pi:C^*(\mathcal{G}) \to B(H)$ be a representation. Suppose that $H_{r(e)}=\oplus_{v\in r(e)}H_v$,  for each $e \in \mathcal{G}^1$. If there exists $n \geq 1$ such that $X_1,\dots,X_n \neq \emptyset$ and $\Big(\bigcup\limits_{e \in \mathcal{G}^1}r(e)\Big) \cup s(\mathcal{G}^1)=\bigcup\limits_{i=1}^{n}(\overline{X_i} \cup I_i)$ then $\pi$ is permutative.
\end{teorema}

\demo
First of all, for $v \in I_0$, take an orthonormal basis $B_v$ of $H_v$.

Secondly, for each $A \in X_1^{\mathrm{fin}}$ and for each $v \in A$, take an orthonormal basis $B_v$ of $H_v$ and let $B_{A}=\bigcup\limits_{v\in A}B_v$. For $A \in X_1^{\mathrm{fin}}$, let $e$ be the corresponding extreme edge associated to $A$ (with $r(e)=A$) and define an orthonormal basis of $H_e$ by $B_e:=\pi(s_e)(B_A)$. For $u \in I_1$, take an orthonormal basis $B_u$ of $H_u$ such that $B_u \supset \bigcup\limits_{f \in s^{-1}(u)}B_f$.

Fix $1 \leq N<n$. Suppose that for each $A\in \bigcup\limits_{i=1}^N X_i^{\mathrm{fin}}$ and for each $v\in \bigcup\limits_{i=1}^N\Big(\Big(\bigcup\limits_{A\in X_i^{\mathrm{fin}}}A\Big)\cup I_i\Big)$ there are orthonormal bases $B_A$ of $H_A$, $B_v$ of $H_v$ such that:

\begin{enumerate}
\item $B_v\supseteq\bigcup\limits_{f\in s^{-1}(v)}B_f$;
\item $B_A =\bigcup\limits_{v\in A}B_v$;
\item $\pi(s_e)(B_A)=B_e$, where $e$ is the extreme edge associated to $A$ with $r(e)=A$.
\end{enumerate}

(Notice that this statement is true when $N=1$, as we shown above).

Fix $A \in X_{N+1}^{\mathrm{fin}}$. For $v \in A$, by Lemma \ref{l1}, we already have a basis $B_f$ of $H_f$, for each $f\in s^{-1}(v)$, and so we take an orthonormal basis $B_v$ of $H_v$ such that $B_v\supset \bigcup\limits_{f \in s^{-1}(v)}B_f$. Let $B_A$ be the basis of $H_A$ defined by $B_A=\bigcup\limits_{v\in A}B_v$. Let $e \in Y_{N+1}$ be the edge associated to $A$, with $r(e)=A$, and define an orthonormal basis of $H_e$ by $B_e:=\pi(s_e)(B_{r(e)})$. Fix $u \in I_{N+1}$. By Lemma \ref{l1}, for each $f\in s^{-1}(u)$, we already have an orthonormal basis $B_f$ and hence we take an orthonormal basis $B_u$ of $H_u$ such that $B_u \supset \bigcup\limits_{f \in s^{-1}(u)}B_f$.

We conclude that for each $A\in \bigcup\limits_{i=1}^{N+1} X_i^{\mathrm{fin}}$, and for each $v\in \bigcup\limits_{i=1}^{N+1}\Big(\Big(\bigcup\limits_{A\in X_i^{\mathrm{fin}}}A\Big)$ $\cup I_i)$, there are orthonormal bases $B_A$ of $H_A$ and $B_v$ of $H_v$ satisfying conditions 1, 2 and 3 above and hence, by induction, for each $A\in \bigcup\limits_{i=1}^{n} X_i^{\mathrm{fin}}$ and $v\in \bigcup\limits_{i=1}^{n}\Big(\Big(\bigcup\limits_{A\in X_i^{\mathrm{fin}}}A\Big)\cup I_i)$,
there are orthonormal bases $B_A$ of $H_A$ and $B_v$ of $H_v$ such that:
\begin{enumerate}\setcounter{enumi}{3}
\item $B_v\supseteq\bigcup\limits_{f\in s^{-1}(v)}B_f$;
\item $B_A=\bigcup\limits_{v\in A} B_v$;
\item $\pi(s_e)(B_A)=B_e$, where $e$ is the extreme edge associated to $A$ with $r(e)=A$.
\end{enumerate}

Fix $e \in \mathcal{G}^1$ such that $r(e) \subset\bigcup\limits_{i=1}^{n}\Big(\Big(\bigcup\limits_{A \in X_i^{\mathrm{fin}}}A\Big) \cup I_i\Big)$ and the orthonormal basis of $H_{e}$ has not been defined yet. Since $H_{r(e)}=\bigoplus\limits_{v\in r(e)}H_v$, define an orthonormal basis of $H_{r(e)}$ by $B_{r(e)}:=\bigcup\limits_{v \in r(e)}B_v$ and define an orthonormal basis of $H_e$ by $B_e:=\pi(s_e)(B_{r(e)})$. For $w \in X_n^{\mathrm{ini}}$, by Lemma~\ref{l1}, we can take an orthonormal basis $B_w$ of $H_w$ such that $B_w \supset \bigcup\limits_{f \in s^{-1}(w)}B_f$.

Now we get an orthonormal basis of the closed subspace $(\oplus_{i=1}^n((\oplus_{A \in X_i^{\mathrm{fin}}}$ $H_A)\oplus (\oplus_{u \in I_i}H_u)))\oplus (\oplus_{w \in X_n^{\mathrm{ini}}}H_w)$ of $H$ satisfying that for each edge $e\in \mathcal{G}^1$, with $r(e)\subseteq \bigcup\limits_{i=1}^n\Big(\Big(\bigcup_{A\in X_i^{\mathrm{fin}}}A\Big)\cup I_i\Big)$, and for each $v\in \bigcup\limits_{i=1}^n\Big(\Big(\bigcup\limits_{A\in X_i^{\mathrm{fin}}}A\Big)\cup I_i\Big)\cup X_n^{\mathrm{ini}}$ there are orthonormal bases $B_{r(e)}$ of $H_{r(e)}$ and $B_v$ of $H_v$ with:

\begin{enumerate}\setcounter{enumi}{6}
\item $B_v\supset\bigcup\limits_{f\in s^{-1}(v)}B_f$;
\item $B_{r(e)}=\bigcup\limits_{v\in r(e)} B_v$;
\item $\pi(s_e)(B_{r(e)})=B_e$.
\end{enumerate}

Fix $1<N\leq n$. Suppose that there exists an orthonormal basis of the closed subspace $\Big(\oplus_{i=1}^{n}((\oplus_{A \in X_i^{\mathrm{fin}}}H_A)\oplus(\oplus_{u \in I_i}H_u))\Big)\oplus\Big(\oplus_{i=N}^{n}\oplus_{w \in X_i^{\mathrm{ini}}}H_w\Big)$ of $H$ satisfying that for each $e\in \mathcal{G}^1$, with $r(e)\subseteq \bigcup\limits_{i=1}^n\Big(\Big(\bigcup\limits_{A\in X_i^{\mathrm{fin}}}A\Big)\cup I_i\Big)\cup\Big(\bigcup\limits_{i=N+1}^n X_i^{\mathrm{ini}}\Big)$, and for each $v\in \bigcup\limits_{i=1}^n\Big(\Big(\bigcup\limits_{A\in X_i^{\mathrm{fin}}}A\Big)\cup I_i\Big)\cup\Big(\bigcup\limits_{i=N}^n X_i^{\mathrm{ini}}\Big)$
there are orthonormal bases $B_{r(e)}$ of $H_{r(e)}$ and $B_v$ of $H_v$ such that conditions 7, 8 and 9 above hold.
%\begin{enumerate}\setcounter{enumi}{12}
%\item $B_v\supset\bigcup\limits_{f\in s^{-1}(v)}B_f$;
%\item $B_{r(e)}=\bigcup\limits_{v\in r(e)}B_v$;
%\item $\pi(s_e)(B_{r(e)})=B_e$;
%\end{enumerate}
(We already showed that the above statement is true when $N=n$, because in this case we consider $\bigcup\limits_{i=N+1}^n X_i^{\mathrm{ini}}$ as being the empty set).

So, fix $e \in \mathcal{G}^1$ such that $r(e) \subset\bigcup\limits_{i=1}^{n}\Big(\Big(\bigcup\limits_{A \in X_i^{\mathrm{fin}}}A\Big) \cup I_i\Big) \cup \bigcup\limits_{i=N}^{n}X_i^{\mathrm{ini}}$ and the orthonormal basis of $H_{e}$ has not been given yet. Since $H_{r(e)}=\bigoplus\limits_{v\in r(e)}H_v$, define an orthonormal basis of $H_{r(e)}$ by $B_{r(e)}:=\bigcup\limits_{v \in r(e)}B_v$ and define an orthonormal basis of $H_e$ by $B_e:=\pi(s_e)(B_{r(e)})$. For $w \in X_{N}^{\mathrm{ini}}$, by Lemma~\ref{l1}, we can take an orthonormal basis $B_w$ of $H_w$ such that $B_w \supset \bigcup\limits_{f \in s^{-1}(w)}B_f$. 

By induction, we get an orthonormal basis of the closed subspace $\Big(\oplus_{i=1}^{n}((\oplus_{A \in X_i^{\mathrm{fin}}}H_A)\oplus(\oplus_{u \in I_i}H_u))\Big) \oplus\Big(\oplus_{i=1}^n\oplus_{w \in X_i^{\mathrm{ini}}}H_w\Big)$ of $H$ satisfying that for each $e\in \mathcal{G}^1$, with $r(e)\subseteq \bigcup\limits_{i=1}^n\Big(\Big(\bigcup\limits_{A\in X_i^{\mathrm{fin}}}A\Big)\cup I_i\Big)\cup\Big(\bigcup\limits_{i=1}^n X_i^{\mathrm{ini}}\Big)=\bigcup\limits_{i=1}^n\overline{X_i}\cup I_i$, (that is, for each $e\in \mathcal{G}$) and for each $v\in \bigcup\limits_{i=1}^n\Big(\Big(\bigcup\limits_{A\in X_i^{\mathrm{fin}}}A\Big)\cup I_i\Big)\cup\Big(\bigcup\limits_{i=1}^n X_i^{\mathrm{ini}}\Big)=\bigcup\limits_{i=1}^n\overline{X_i}\cup I_i$ there are orthonormal bases $B_{r(e)}$ of $H_{r(e)}$ and $B_v$ of $H_v$ such that:

\begin{enumerate}\setcounter{enumi}{9}
\item $B_v\supset\bigcup\limits_{f\in s^{-1}(v)}B_f$;
\item $B_{r(e)}= \bigcup\limits_{v\in r(e)} B_v$;
\item $\pi(s_e)(B_{r(e)})=B_e$, for all $e \in \mathcal{G}^1$;
\end{enumerate}

Finally, let $B$ an orthonormal basis of $H$ such that $B_v\subseteq B$ for each $v\in G^0$. Then $B$ satisfies Conditions~(1)--(4) of Definition~\ref{define permutative rep} and therefore we are done. \fim

%================================================================================================================================
\section{An Interlude: General Cuntz-Krieger Theorem for Ultragraph C*-algebras}
%================================================================================================================================

In this section we make a pause on the theory of branching systems to develop a generalized Cuntz-Krieger theorem for ultragraph C*-algebras. This result is of independent interest and key to the development of the theory of faithful representations of ultragraph C*-algebras.

\begin{lema}\label{I_SH=I_SH1+I_SH2}
Let $\mathcal{G}$ be an ultragraph and $\mathcal{SH}_1, \mathcal{SH}_2$ be hereditary and saturated subsets of $\mathcal{G}^0$. Define $X:=\{A\cup B:A \in \mathcal{SH}_1,B\in \mathcal{SH}_2\}$ ($X$ is a hereditary subset of $\mathcal{G}^0$). Let $\mathcal{SH}$ be the smallest hereditary and saturated subset of $\mathcal{G}^0$ containing $X$. Then $I_{\mathcal{SH}}=I_{\mathcal{SH}_1}+I_{\mathcal{SH}_2}$. Moreover, if $\mathcal{SH}_1\cap\mathcal{SH}_2=\emptyset$ then $I_{\mathcal{SH}_1}I_{\mathcal{SH}_2}=0$ (Recall that $I_{\mathcal{SH}}$ was defined in Lemma \ref{define ideal I_H}).
\end{lema}

\demo
Let $X_0:=X$. For $n \geq 0$, denote by $S_n, X_{n+1}$ the sets defined as in Lemma~\ref{smallest here and sat} such that $\mathcal{SH}=\bigcup_{n=0}^{\infty}X_n$.

It is straightforward to see that $I_{\mathcal{SH}_1}+I_{\mathcal{SH}_2} \subset I_{\mathcal{SH}}$. We show that $I_{\mathcal{SH}}\subset I_{\mathcal{SH}_1}+I_{\mathcal{SH}_2}$. By the definition of $I_{\mathcal{SH}}$, we only need to prove that $p_A \in I_{\mathcal{SH}_1}+I_{\mathcal{SH}_2}$ for all $A \in \mathcal{SH}$. For $A \in \mathcal{SH}_1, B \in \mathcal{SH}_2$ ($A \cup B \in X_0$), we have that $A \cap B \in \mathcal{SH}_2$ because $\mathcal{SH}_2$ is hereditary. So $p_{A \cup B}=p_A+p_B-p_{A \cap B} \in I_{\mathcal{SH}_1}+I_{\mathcal{SH}_2}$.

Suppose that for $N \geq 0$ and for $A \in \bigcup_{n=0}^{N}X_n$, we have that $p_A \in I_{\mathcal{SH}_1}+I_{\mathcal{SH}_2}$. Fix $A \in X_N, F \in S_N$ ($A \cup F \in X_{N+1}$). We have $p_{A \cup F}=p_A+p_F-p_Ap_F$. By hypothesis $p_A \in I_{\mathcal{SH}_1}+I_{\mathcal{SH}_2}$. For $v \in F$, we have $p_v=\sum_{e \in s^{-1}(v)}s_es_e^*=\sum_{e \in s^{-1}(v)}s_ep_{r(e)}s_e^* \in I_{\mathcal{SH}_1}+I_{\mathcal{SH}_2}$ because $0<\vert s^{-1}(v)\vert<\infty$ and $r(s^{-1}(v)) \subset X_N$. Since $p_F=\sum_{v \in F}p_v$, we obtain that $p_F \in I_{\mathcal{SH}_1}+I_{\mathcal{SH}_2}$. So $p_{A \cup F} \in I_{\mathcal{SH}_1}+I_{\mathcal{SH}_2}$ and hence, for $A \in \bigcup_{n=0}^{N+1}X_n$, we have that $p_A \in I_{\mathcal{SH}_1}+I_{\mathcal{SH}_2}$ and the induction step is verified. We deduce that for all $A \in \mathcal{SH}, p_A \in I_{\mathcal{SH}_1}+I_{\mathcal{SH}_2}$ and therefore $I_{\mathcal{SH}}=I_{\mathcal{SH}_1}+I_{\mathcal{SH}_2}$.

Now suppose that $\mathcal{SH}_1\cap\mathcal{SH}_2=\emptyset$. By Lemma~\ref{define ideal I_H}, we only need to show that $p_As_\beta^*s_\gamma p_B=0$, for all $A \in \mathcal{SH}_1,B \in \mathcal{SH}_2, \beta,\gamma \in \mathcal{G}^*$, which follows from a straightforward calculation. \fim

\begin{lema}\label{gauge-inv ideal structure}
Let $\mathcal{G}$ be an ultragraph and $I$ be a nonzero, gauge-invariant, closed, two-sided ideal of $C^*(\mathcal{G})$. Then there exists $v \in G^0$ such that $p_v \in I$.
\end{lema}

\demo
If $I=C^*(\mathcal{G})$ then the statement is trivially valid. Suppose that $I \neq C^*(\mathcal{G})$. Denote by $\{p_A,s_e:A\in\mathcal{G}^0,e\in \mathcal{G}^1\}$ the Cuntz-Krieger $\mathcal{G}$-family generating $C^*(\mathcal{G})$. Suppose that $p_A \notin I$ for all nonempty $A \in \mathcal{G}^0$, for a contradiction. Consider a family of partial isometries $\{s_e+I:e \in \mathcal{G}^1 \}$ and a family of projections $\{p_A+I: A \in \mathcal{G}^0\}$ in the C*-algebra $C^*(\mathcal{G})/I$. It is straightforward to see that $\{p_A+I,s_e+I:A\in\mathcal{G}^0,e\in \mathcal{G}^1\}$ is a Cuntz-Krieger $\mathcal{G}$-family. By the universal property of $C^*(\mathcal{G})$, there exists a homomorphism $\pi:C^*(\mathcal{G}) \to C^*(\mathcal{G})/I$ such that $\pi(p_A)=p_A+I$ and $\pi(s_e)=s_e+I$, for all $A \in \mathcal{G}^0, e \in \mathcal{G}^1$. Notice that $\pi$ is indeed the quotient map. Since $I \neq C^*(\mathcal{G})$ and $p_A \notin I$ for all nonempty $A \in \mathcal{G}^0$, we have $p_A+I\neq 0$ for all nonempty $A \in \mathcal{G}^0$. Since $I$ is gauge-invariant, the gauge action on $C^*(\mathcal{G})$ induces a gauge action on $C^*(\mathcal{G})/I$. So \cite[Theorem~6.8]{Tomforde:JOT03} implies that $\pi$ is an isomorphism. However this is impossible because $I$ is nonzero. Therefore we are done. \fim

The following lemma is a generalization of \cite[Lemma~1.1]{Szyma'nski:IJM02}. For the reader convenience, we recall that a closed ideal $I$ in a C*-algebra $A$ is called essential if it has non-zero intersection with every other non-zero closed ideal, or, equivalently, if $aI=0$ implies $a=0$ for all $a\in A$.

\begin{lema}\label{criterion of ideal being essential}
Let $\mathcal{G}$ be an ultragraph and let $\mathcal{SH}$ be a hereditary and saturated subset of $\mathcal{G}^0$. Then the closed two-sided ideal $I_{\mathcal{SH}}$ is essential if and only if for $v \in G^0 \setminus\Big(\bigcup_{A \in \mathcal{SH}}A\Big)$, there exists $\alpha \in \mathcal{G}^*\setminus\mathcal{G}^0$, such that $s(\alpha_1)=v, r(\alpha_{\vert\alpha\vert}) \cap\Big(\bigcup_{A \in \mathcal{SH}}A\Big)\neq\emptyset$.
\end{lema}

\demo
Suppose that $I_{\mathcal{SH}}$ is essential. Fix $v \in G^0 \setminus\Big(\bigcup_{A \in \mathcal{SH}}A\Big)$. Then $p_vI_{\mathcal{SH}}\neq 0$. By Lemma~\ref{define ideal I_H}, there exist $\alpha,\beta\in\mathcal{G}^*$, $A \in \mathcal{SH}$ such that $p_vs_\alpha p_As_\beta^*\neq0$. A straightforward calculation gives $\alpha \in\mathcal{G}^*\setminus\mathcal{G}^0$ such that $s(\alpha_1)=v, r(\alpha_{\vert\alpha\vert}) \cap A\neq\emptyset$.

Conversely, suppose that for any $v \in G^0 \setminus\Big(\bigcup_{A \in \mathcal{SH}}A\Big)$, there exists $\alpha \in \mathcal{G}^*\setminus\mathcal{G}^0$, such that $s(\alpha_1)=v, r(\alpha_{\vert\alpha\vert}) \cap\Big(\bigcup_{A \in \mathcal{SH}}A\Big)\neq\emptyset$. Suppose that $I_{\mathcal{SH}}$ is not essential, for a contradiction. Then $I_{\mathcal{SH}}^{\perp}\neq 0$. Since $I_{\mathcal{SH}}$ is gauge-invariant by Lemma~\ref{define ideal I_H}, $I_{\mathcal{SH}}^\perp$ is gauge-invariant as well. By Lemma~\ref{gauge-inv ideal structure}, there exists $v \in G^0$ such that $p_{\{v\}} \in I_{\mathcal{SH}}^\perp$. Notice that $v \in  G^0 \setminus\Big(\bigcup_{A \in \mathcal{SH}}A\Big)$. So there exists $\alpha \in \mathcal{G}^*\setminus\mathcal{G}^0$, such that $s(\alpha_1)=v, r(\alpha_{\vert\alpha\vert}) \cap\Big(\bigcup_{A \in \mathcal{SH}}A\Big)\neq\emptyset$. Take $w \in r(\alpha_{\vert\alpha\vert}) \cap\Big(\bigcup_{A \in \mathcal{SH}}A\Big)$. We have $p_{\{v\}} s_\alpha p_{r(\alpha)}p_{\{w\}}=0$ because $p_{\{v\}} \in I_{\mathcal{SH}}^\perp, p_{\{w\}} \in I_{\mathcal{SH}}$. However, $(p_{\{v\}} s_\alpha p_{r(\alpha)}p_{\{w\}} )^*p_{\{v\}} s_\alpha p_{r(\alpha)}p_{\{w\}} =p_{\{w\}} \neq 0$, which is a contradiction. Hence $I_{\mathcal{SH}}$ is essential. \fim

The following theorem is a generalization of \cite[Theorem~1.2]{Szyma'nski:IJM02}.

\begin{teorema}\label{general Cuntz-Krieger for ultragraph}
Let $\mathcal{G}$ be an ultragraph, $\mathcal{A}$ be a C*-algebra and let $\varphi:C^*(\mathcal{G}) \to \mathcal{A}$ be a homomorphism. Then $\varphi$ is injective if and only if the following hold.
\begin{enumerate}
\item $\varphi(p_A) \neq 0$ for all nonempty $A \in \mathcal{G}^0$;
\item for any simple cycle $\alpha$, without exits, the spectrum of $\varphi(s_\alpha)$ contains the unit circle.
\end{enumerate}
\end{teorema}

\demo
First of all, suppose that $\varphi$ is injective. Then Condition~(1) is trivially true. We check Condition~(2). Fix a simple cycle $\alpha$ without exits. Then \cite[Lemma~3.8]{Tomforde:IUMJ03} yields that the spectrum of $s_\alpha$ contains the unit circle and, since $\varphi$ is injective, \cite[Corollary~II.1.6.7]{Blackadar:Operatoralgebras06} implies that the spectrum of $\varphi(s_\alpha)$ contains the unit circle.

Conversely, suppose that Conditions~(1), (2) hold. Let $X_1:=\{s(\alpha):\alpha \text{ is a simple cycle without exits in } \mathcal{G} \}$. Denote by $\mathcal{H}_1$ the class of all finite subsets of $X_1$, which is a hereditary subset of $\mathcal{G}^0$. Denote by $\mathcal{SH}_1$ the smallest hereditary and saturated subset of $\mathcal{G}^0$ containing $\mathcal{H}_1$ (see Lemma~\ref{smallest here and sat}). Define
\begin{align*}
\mathcal{SH}_2:=\{A\in \mathcal{G}^0: \text{for any } \alpha\in\mathcal{G}^*\setminus\mathcal{G}^0,  s(\alpha)\in A &\implies
\\&r(\alpha) \cap \Big(\bigcup_{B\in\mathcal{SH}_1}B\Big)=\emptyset\},
\end{align*}
which is a hereditary and saturated subset of $\mathcal{G}^0$. Define $X:=\{A\cup B:A \in \mathcal{SH}_1,B\in \mathcal{SH}_2\}$ which is a hereditary subset of $\mathcal{G}^0$. Denote by $\mathcal{SH}$ the smallest hereditary and saturated subset of $\mathcal{G}^0$ containing $X$.

By Lemma~\ref{I_SH=I_SH1+I_SH2}, $I_{\mathcal{SH}}=I_{\mathcal{SH}_1}+I_{\mathcal{SH}_2}$. Notice that $\mathcal{SH}_1 \cap \mathcal{SH}_2=\emptyset$ and hence Lemma~\ref{I_SH=I_SH1+I_SH2} implies that $I_{\mathcal{SH}_1}I_{\mathcal{SH}_2}=0$.

Fix $v \in G^0 \setminus\Big(\bigcup_{A \in \mathcal{SH}}A\Big)$. Then $s^{-1}(v) \neq \emptyset$ since otherwise $\{v\} \in \mathcal{SH}_2$. We deduce that there exists $\alpha \in \mathcal{G}^*\setminus\mathcal{G}^0$, such that $s(\alpha_1)=v, r(\alpha_{\vert\alpha\vert}) \cap\Big(\bigcup_{A \in \mathcal{SH}}A\Big)\neq\emptyset$, because otherwise $v \in \mathcal{SH}_2$. By Lemma~\ref{criterion of ideal being essential}, $I_{\mathcal{SH}}$ is essential. Hence, in order to show that $\varphi$ is injective, we only need to show that $\ker(\varphi) \cap I_{\mathcal{SH}}=0$. Since $I_{\mathcal{SH}}=I_{\mathcal{SH}_1}+I_{\mathcal{SH}_2}$ and $I_{\mathcal{SH}_1}I_{\mathcal{SH}_2}=0$, it is equivalent to show that $\ker(\varphi) \cap I_{\mathcal{SH}_1}=0$ and $\ker(\varphi) \cap I_{\mathcal{SH}_2}=0$.

Now we prove that $\ker(\varphi) \cap I_{\mathcal{SH}_1}=0$. By Lemma \ref{define ideal I_H},
\[
I_{\mathcal{SH}_1}=\overline{\mathrm{span}}\{s_\alpha p_A s_\beta^*:\alpha,\beta \in \mathcal{G}^*,A \in \mathcal{SH}_1\}.
\]
Denote by $\mathcal{B}$ the C*-subalgebra of $I_{\mathcal{SH}_1}$ generated by
\begin{align*}
\Big\{p_{\{s(\alpha_1),\dots,s(\alpha_{\vert\alpha\vert})\}}I_{\mathcal{SH}_1}p_{\{s(\alpha_1),\dots,s(\alpha_{\vert\alpha\vert})\}}:\alpha \text{ is a simple cycle without exits } \Big\}.
\end{align*}
Denote a closed subspace of $I_{\mathcal{SH}_1}$ by
\[
Y:=\overline{\mathrm{span}}\{s_\alpha p_A s_\beta^*:s(\alpha) \subset X_1, A \in \mathcal{SH}_1\}.
\]
We show that $YY^* \subset \mathcal{B}$. Notice that each element in the spanning set of $YY^*$ has the form $s_\alpha p_A s_{\alpha'}^*$, where $s(\alpha),s(\alpha') \subset X_1$ and $A \in \mathcal{SH}_1$.

Case $1$. $\vert\alpha\vert,\vert\alpha'\vert\geq 1$. Since $s(\alpha), s(\alpha') \subset X_1$, we deduce that $r(\alpha), r(\alpha')$ are both singletons in $X_1$. We may assume that $r(\alpha)=r(\alpha')$ because otherwise $s_\alpha p_A s_{\alpha'}^*=0$. So $s(\alpha),s(\alpha')$ are on the same cycle without exits, say $\beta$. Hence
\[
s_\alpha p_A s_{\alpha'}^*=p_{\{s(\beta_1),\dots,s(\beta_{\vert\beta\vert})\}}s_\alpha p_A s_{\alpha'}^* p_{\{s(\beta_1),\dots,s(\beta_{\vert\beta\vert})\}} \in \mathcal{B}.
\]

Case $2$. $\vert\alpha\vert=0, \vert\alpha'\vert\geq 1$. Since $s(\alpha') \subset X_1, r(\alpha')$ is a singleton and $r(\alpha'),s(\alpha')$ are on the same cycle without exits, say $\beta$. We may assume that $\alpha \cap A \cap r(\alpha') \neq 0$, because otherwise $s_\alpha p_A s_{\alpha'}^*=0$. So
\[
s_\alpha p_A s_{\alpha'}^*=s_{\alpha'}^*=p_{\{s(\beta_1),\dots,s(\beta_{\vert\beta\vert})\}}s_\alpha p_A s_{\alpha'}^* p_{\{s(\beta_1),\dots,s(\beta_{\vert\beta\vert})\}} \in \mathcal{B}.
\]

Case $3$. $\vert\alpha\vert\geq 1, \vert\alpha'\vert =0$. It follows from Case $2$.

Case $4$. $\vert\alpha\vert, \vert\alpha'\vert =0$. Since $\mathcal{H}_1$ consists of all finite subsets of $X_1$, Lemma~\ref{smallest here and sat} yields that $\alpha \cap A \cap \alpha'$ is a finite subset of $X_1$. So it is easy to see that $s_\alpha p_A s_{\alpha'}^*\in \mathcal{B}$.

It is straightforward to check that $\mathcal{B}Y,YI_{\mathcal{SH}_1} \subset Y; Y^*Y \subset I_{\mathcal{SH}_1}$; and $\overline{\mathrm{span}}YY^*=\mathcal{B},\overline{\mathrm{span}}Y^*Y=I_{\mathcal{SH}_1}$. So by Lemma~\ref{imprimitivity bimodule}, in order to prove that $\ker(\varphi) \cap I_{\mathcal{SH}_1}=0$, it is sufficient to show that $\ker(\varphi) \cap \mathcal{B}=0$. By the construction of $\mathcal{B}$, we only have to show that $\ker(\varphi) \cap (p_{\{s(\alpha_1),\dots,s(\alpha_{\vert\alpha\vert})\}}I_{\mathcal{SH}_1}$ $p_{\{s(\alpha_1),\dots,s(\alpha_{\vert\alpha\vert})\}})=0$ for each cycle without exits $\alpha$. Fix a cycle without exits $\alpha$. Let $\Lambda$ be the smallest hereditary and saturated subset of $\mathcal{G}^0$ containing $\{\{s_{\alpha_i}\}\}_{i=1}^{\vert\alpha\vert}$. Notice that
\[
p_{\{s(\alpha_1),\dots,s(\alpha_{\vert\alpha\vert})\}}I_{\mathcal{SH}_1}p_{\{s(\alpha_1),\dots,s(\alpha_{\vert\alpha\vert})\}}=p_{\{s(\alpha_1),\dots,s(\alpha_{\vert\alpha\vert})\}}I_{\Lambda}p_{\{s(\alpha_1),\dots,s(\alpha_{\vert\alpha\vert})\}}.
\]
By \cite[Lemma~3.8]{Tomforde:IUMJ03}, $p_{\{s(\alpha_1),\dots,s(\alpha_{\vert\alpha\vert})\}}I_{\mathcal{SH}_1}p_{\{s(\alpha_1),\dots,s(\alpha_{\vert\alpha\vert})\}}$ is isomorphic with $M_{\vert\alpha\vert}(\mathbb{C}) \otimes C(\mathbb{T})$. Now, \cite[Exercises~1.11.42]{Lin:introductiontoclassification01} yields that 
\[
\ker(\varphi) \cap(p_{\{s(\alpha_1),\dots,s(\alpha_{\vert\alpha\vert})\}}I_{\mathcal{SH}_1} p_{\{s(\alpha_1),\dots,s(\alpha_{\vert\alpha\vert})\}})=0.
\]
So $\ker(\varphi) \cap I_{\mathcal{SH}_1}=0$.

Finally we prove that $\ker(\varphi) \cap I_{\mathcal{SH}_2}=0$. By \cite[Lemma~3.5]{Tomforde:IUMJ03},
\[
I_{\mathcal{SH}_2}=\overline{\mathrm{span}}\{s_\alpha p_A s_\beta^*:\alpha,\beta \in \mathcal{G}^*,A \in \mathcal{SH}_2\}.
\]
Denote by $\mathcal{C}$ the C*-subalgebra of $I_{\mathcal{SH}_2}$ generated by
\[
\Big\{p_A, s_e:A \in \mathcal{SH}_2, e \in s^{-1}\Big(\bigcup_{A \in \mathcal{SH}_2}A\Big) \Big\}.
\]
Denote a closed subspace of $I_{\mathcal{SH}_2}$ by
\[
Z:=\overline{\mathrm{span}}\{s_\alpha p_A s_\beta^*:s(\alpha), A \in \mathcal{SH}_2\}.
\]
It is straightforward to verify that $\mathcal{C}Z,ZI_{\mathcal{SH}_2} \subset Z; ZZ^* \subset \mathcal{C},Z^*Z \subset I_{\mathcal{SH}_2}$; and $\overline{\mathrm{span}}ZZ^*=\mathcal{C},\overline{\mathrm{span}}Z^*Z=I_{\mathcal{SH}_2}$. By Lemma~\ref{imprimitivity bimodule}, in order to prove that $\ker(\varphi) \cap I_{\mathcal{SH}_2}=0$, it is enough to show that $\ker(\varphi) \cap \mathcal{C}=0$.

Define an ultragraph 
\begin{align*}
\mathcal{F}&=(F^0,\mathcal{F}^1,r_{\mathcal{F}},s_{\mathcal{F}})
\\&:=(G^0,\mathcal{G}^1 \setminus \{\alpha_1,\dots,\alpha_{\vert\alpha\vert}:\alpha \text{ is a cycle without exits}\},r,s).
\end{align*}
Denote by $\Big\{q_A, t_e:A \in \mathcal{F}^0,e\in \mathcal{F}^1 \Big\}$ the Cuntz-Krieger $\mathcal{F}$-family generating $C^*(\mathcal{F})$. Denote by $\mathcal{D}$ the C*-subalgebra of $C^*(\mathcal{G})$ generated by
\[
\Big\{p_A, s_e:A \in \mathcal{F}^0,e\in \mathcal{F}^1 \Big\}. 
\]
Lemma~\ref{concrete des of G^0} gives $\mathcal{C} \subset \mathcal{D}$. It is straightforward to verify that $\Big\{p_A, s_e:A \in \mathcal{F}^0,e\in \mathcal{F}^1 \Big\}$ is a Cuntz-Krieger $\mathcal{F}$-family. Let $\pi:C^*(\mathcal{F}) \to \mathcal{D}$ be the surjective homomorphism induced from the universal property of $C^*(\mathcal{F})$. By Condition~(1), we have $\varphi\circ\pi(q_A)\neq 0$ for all nonempty $A \in \mathcal{F}^0$. Notice that $\mathcal{F}$ satisfies Condition~(L). By Theorem~\ref{Cuntz-Krieger uni thm for ultragraph alg}, the homomorphism $\varphi \circ \pi$ is injective. So $\ker(\varphi) \cap \mathcal{D}=0$. Since $\mathcal{C} \subset \mathcal{D}, \ker(\varphi) \cap \mathcal{C}=0$. Hence $\ker(\varphi) \cap I_{\mathcal{SH}_2}=0$ and therefore we are done. \fim

\begin{rmk}
When $\mathcal{G}$ is a directed graph, as shown in \cite[Theorem~1.2]{Szyma'nski:IJM02}, the C*-subalgebra $\mathcal{D}$ was not needed because $\mathcal{C}$ is automatically isomorphic to a graph algebra whose underlying graph satisfies Condition~(L). More precisely, one could construct a directed graph satisfying Condition~(L) by letting $\mathcal{F}':=(\bigcup_{A \in \mathcal{SH}_2}A,s^{-1}\Big(\bigcup_{A \in \mathcal{SH}_2}A\Big),r,s)$. Then the generator of $\mathcal{C}$ is a Cuntz-Krieger $\mathcal{F}'$-family, so there exists a surjective homomorphism $\pi':C^*(\mathcal{F}') \to \mathcal{C}$. Hence $\varphi\circ\pi'$ is injective, which yields that $\ker(\varphi) \cap C=0$. However, this argument will not work when dealing with ultragraphs in general.

For example, consider the ultragraph $\mathcal{G}=(\{u,v,w_n,x,y\}_{n=1}^{\infty}, \{e,f,g,h\}$, $r,s)$, where $r(e)=\{w_n,x\}_{n=1}^{\infty};s(e)=u;r(f)=\{w_n,y\}_{n=1}^{\infty};s(f)=v;r(g)=s(g)=x;r(h)=s(h)=y$. In this case, $\mathcal{SH}_2=\{ \{w_n\}_{n=1}^{\infty}, \text{ all finite subsets }$ $\text{of } \{w_n\}_{n=1}^{\infty}\}$, but $s^{-1}(\bigcup_{A \in \mathcal{SH}_2}A)=\emptyset$.
\end{rmk}

%================================================================================================================================
\section{Faithful Representations of Ultragraph C*-algebras via Branching Systems}
%================================================================================================================================

To finish the paper we present some results regarding faithfulness of representations of ultragraph C*-algebras arising from branching systems. 

\begin{proposicao}\label{faithfulrep}
Let $\mathcal{G}$ be an ultragraph satisfying Condition~(L). Fix a $\mathcal{G}$-branching system $\{R_e,D_A,f_e\}_{e\in \mathcal{G}^1, A\in \mathcal{G}^0}$ on a measure space $(X,\mu)$. Suppose that $\mu(D_A) \neq 0$ for any nonempty set $A \in \mathcal{G}^0$. Then the representation $\pi:C^*(\mathcal{G}) \to B(\mathcal{L}^2(X,\mu))$ induced from the branching system (see Theorem~\ref{repinducedbybranchingsystems}) is faithful.
\end{proposicao}

\demo
Fix $A \in \mathcal{G}^0$ with $A \neq \emptyset$. Take an arbitrary $v \in A$ and a measurable set $S \subset D_v$ such that $0<\mu(S)<\infty$. Since $\pi(p_A)(\chi_{S})=\chi_{D_A}\chi_S=\chi_S \neq 0$, we get that $\pi(p_A) \neq 0$ and hence Theorem~\ref{Cuntz-Krieger uni thm for ultragraph alg} implies that $\pi$ is faithful. \fim

The following theorem is the converse of the Cuntz-Krieger uniqueness theorem for ultragraph C*-algebras, which is a generalization of \cite[Theorem~3.5]{GLR}.

\begin{teorema}
Let $\mathcal{G}$ be an ultragraph which does not satisfy Condition~(L). Then there exists a $\mathcal{G}$-branching system $\{R_e,D_A,f_e\}_{e\in \mathcal{G}^1, A\in \mathcal{G}^0}$ on a measure space $(X,\mu)$ with the induced representation $\pi:C^*(\mathcal{G}) \to B(\mathcal{L}^2(X,\mu))$ (see Theorem~\ref{repinducedbybranchingsystems}) such that $\pi(p_A)\neq 0$, for all nonempty $A \in \mathcal{G}^0$, and $\pi$ is not faithful.
\end{teorema}

\demo
Since $\mathcal{G}$ does not satisfy Condition~(L), there exists a cycle $\alpha=(\alpha_1,\dots,\alpha_n)$ such that $\vert r(\alpha_i) \vert=1, s^{-1}(s(\alpha_i))=\{\alpha_i\}$ for all $i=1,\dots,n$, and $\alpha_i \neq \alpha_j$ if $i \neq j$. We enumerate the edge set as $\mathcal{G}^1=\{\alpha_1,\dots,\alpha_n,e_{n+1},\dots\}$, and enumerate the vertex set as $G^0=\{s(\alpha_1),\dots,s(\alpha_n),v_{n+1},\dots\}$. By Theorem~\ref{existenceofabranchingsystem}, there is a $\mathcal{G}$-branching system on $(\mathbb{R},\mu)$ denoted by $\{R_e,D_A,f_e:e\in \mathcal{G}^1, A\in \mathcal{G}^0\}$, where $\mu$ is the Lebesgue measure on all Borel sets of $\mathbb{R}$, such that for each $i=1,\dots,n$, we have $D_{s(\alpha_i)}=R_{\alpha_i}=[i-1,i]$ and $f_{\alpha_i}$ is the increasing bijective linear map. Notice that $f_\alpha=\id$ and so $\Phi_{f_\alpha}\equiv 1$ on $[0,1]$. So $\pi(s_\alpha^*)=\pi(p_{s(\alpha_1)})$. By Theorem~\ref{existenceofabranchingsystem}, we deduce that $\pi(p_A)\neq 0$ for all nonempty $A \in \mathcal{G}^0$, and then in particular, $\pi(p_{s(\alpha_1)})\neq 0$.

Suppose that $\pi$ is faithful, for a contradiction. By the universal property of $C^*(\mathcal{G})$, there exists a gauge action $\gamma$ on $\pi(C^*(\mathcal{G}))$. So for each $z \in \mathbb{T}$ we have that
\[
\pi(p_{s(\alpha_1)})=\gamma_z(\pi(p_{s_{\alpha_1}}))=\gamma_z(\pi(s_\alpha^*))=\overline{z}^n\pi(s_\alpha^*)=\overline{z}^n\pi(p_{s(\alpha_1)}),
\]
which is impossible, since $\pi(p_{s(\alpha_1)})\neq 0$. Therefore $\pi$ is not faithful.
\fim

\begin{teorema}\label{a criterion of faithful rep}
Let $\mathcal{G}$ be an ultragraph, let $\{R_e,D_A,f_e:e\in \mathcal{G}^1, A\in \mathcal{G}^0\}$ be a $\mathcal{G}$-branching system on a measure space $(X,\mu)$ such that $\mu(D_A) \neq 0$, for all nonempty $A \in \mathcal{G}^0$, and let $\pi:C^*(\mathcal{G}) \to B(\mathcal{L}^2(X,\mu))$ be the representation induced from the branching system (see Theorem~\ref{repinducedbybranchingsystems}). Suppose that for any simple cycle $\alpha=(\alpha_i)_{i=1}^{n}$ without exits, and for any finite subset $\mathcal{F}$ of $\mathbb{N}$, there exists a measurable subset $E$ of $D_{s(\alpha_1)}$, with $\mu(E) \neq 0$, such that $f_{n\alpha}(E) \cap E\stackrel{\mu-a.e.}{=}\emptyset$ for all $n \in \mathcal{F}$ ($n\alpha$ stands for $\overbrace{\alpha\dots\alpha}^n$). Then $\pi$ is faithful.
\end{teorema}

\demo
For each nonempty set $A \in \mathcal{G}^0$, since $\mu(D_A) \neq 0$, we have that $\pi(p_A) \neq 0$. Fix a simple cycle $\alpha=(\alpha_i)_{i=1}^{n}$ without exits, and fix a finite subset $\mathcal{F}$ of $\mathbb{N}$. By Theorem~\ref{general Cuntz-Krieger for ultragraph}, we only need to show that the spectrum of $\pi(s_\alpha)$ contains the unit circle. Since $\alpha$ is a simple cycle without exits, we have that $\pi(s_\alpha)^*\pi(s_\alpha)=\pi(s_\alpha)\pi(s_\alpha)^*=\pi(p_{s(\alpha_1)})$. So there exists a unique homomorphism $h:C(\mathbb{T}) \to C^*(\pi(s_{\alpha}))$, such that $h(I)=\pi(p_{s(\alpha_1)})$ and $h(u)=\pi(s_\alpha)$, where $u$ is the unitary generating $C(\mathbb{T})$. Since the spectrum of $u$ is the unit circle, by \cite[Corollary~II.1.6.7]{Blackadar:Operatoralgebras06}, to prove that the spectrum of $\pi(s_\alpha)$ in $\pi(C^*(\mathcal{G}))$ contains the unit circle, it is sufficient to prove that $h$ is an isomorphism. Notice that there exists a conditional expectation $\Phi:C(\mathbb{T}) \to C(\mathbb{T})$ such that $\Phi(u^n(u^*)^m)=\delta_{n,m}I$. By \cite[Proposition~3.11]{Katsura:CJM03}, to show that $h$ is faithful, we only need to construct a conditional expectation $\Psi:C^*(\pi(s_{\alpha})) \to C^*(\pi(s_{\alpha}))$ such that $\Psi(\pi(s_{n\alpha})\pi(s_{m\alpha})^*)=\delta_{n,m}\pi(p_{s(\alpha_1)})$. It is sufficient to show that
\[
\vert z\vert\leq\Big\Vert z\pi(p_{s(\alpha_1)})+\sum_{n \in \mathcal{F}}z_n\pi(s_{n\alpha})+\sum_{m \in \mathcal{F}'}z_m'\pi(s_{m\alpha})^*\Big\Vert
\]
where $\mathcal{F},\mathcal{F}'$ are finite subsets of $\mathbb{N}$.

Fix finite subsets $\mathcal{F},\mathcal{F}' \subset \mathbb{N}$. By the assumption of the theorem, there exists a measurable subset $E$ of $D_{s(\alpha_1)}$ with $\mu(E) \neq 0$, such that $f_{n\alpha}(E) \cap E\stackrel{\mu-a.e.}{=}\emptyset$ for all $n \in \mathcal{F} \cup \mathcal{F}'$. Take an arbitrary function $\phi \in \mathcal{L}^2(X,\mu)$ with $\Vert\phi\Vert = 1$ and $\supp(\phi) \stackrel{\mu-a.e.}{\subset}E$. Then $\pi(s_{n\alpha})(\phi)(x)=0$ and $\pi(s_{m\alpha})^*(\phi)(x)=0$ for almost every $x\in E$ and for all $n\in \mathcal{F}$ and $m\in \mathcal{F}'$.
Then
\begin{align*}
\Big\Vert z\pi(p_{s(\alpha_1)})(\phi)&+\sum_{n \in \mathcal{F}}z_n\pi(s_{n\alpha})(\phi)+\sum_{m \in \mathcal{F}'}z_m'\pi(s_{m\alpha})^*(\phi)\Big\Vert^2
\\&\geq\int_E  \Big\vert z\pi(p_{s(\alpha_1)})(\phi)+\sum_{n \in \mathcal{F}}z_n\pi(s_{n\alpha})(\phi)+\sum_{m \in \mathcal{F}'}z_m'\pi(s_{m\alpha})^*(\phi) \Big\vert^2 \
\\& = \int_E \Big\vert z\pi(p_{s(\alpha_1)})(\phi)\Big\vert^2\,\mathrm{d}\mu
\\&=\vert z\vert^2.
\end{align*}
So $\vert z\vert\leq\Big\Vert z\pi(p_{s(\alpha_1)})+\sum\limits_{n \in \mathcal{F}}z_n\pi(s_{n\alpha})+\sum\limits_{m \in \mathcal{F}'}z_m'\pi(s_{m\alpha})^*\Big\Vert$ and we are done. \fim

\section*{Acknowledgments}

The second author would like to thank the rest of the authors for the continuation of the collaboration and for their extremely friendly and patient correspondences. Finally the second author appreciates the support from Research Center for Operator Algebras of East China Normal University.

\vspace{1.5pc}
\begin{center}
Daniel Gon\c{c}alves (daemig@gmail.com) and Danilo Royer (danilo.royer@ufsc.br)\\

Departamento de Matem\'{a}tica - Universidade Federal de Santa Catarina, Florian\'{o}polis, 88040-900, Brazil

\vspace{1.5pc}

Hui Li (lihui8605@hotmail.com)\\

Research Center for Operator Algebras, Department of Mathematics, East China Normal University, 3663 Zhongshan North Road, Putuo District, Shanghai 200062, China
\end{center}


\begin{thebibliography}{15}
\bibitem{BatesHongEtAl:IJM02} T. Bates, J.H. Hong, I. Raeburn, and W. Szyma{\'n}ski, \emph{The ideal structure of the {$C\sp *$}-algebras of infinite graphs}, Illinois J. Math. \textbf{46} (2002), 1159--1176.
\bibitem{Blackadar:Operatoralgebras06} B. Blackadar, Operator algebras, Theory of $C\sp *$-algebras and von Neumann algebras, Operator Algebras and Non-commutative Geometry, III, Springer-Verlag, Berlin, 2006, xx+517.
\bibitem{MR1465320} O. Bratteli and P.E.T. Jorgensen, \emph{Isometries, shifts, {C}untz algebras and multiresolution wavelet analysis of scale {$N$}}, Integral Equations Operator Theory \textbf{28} (1997), 382--443.
\bibitem{BJ} O. Bratteli and P.E.T. Jorgensen, \emph{Iterated function Systems and Permutation Representations of the Cuntz algebra}, Memories of the American Mathematical Society, \textbf{139},
(1999).
\bibitem{BNS} J.H. Brown, G. Nagy and S. Reznikoff, \emph{Generalized {C}untz-{K}rieger uniqueness theorem for higher-rank graphs}, Journal of Functional Analysis, \textbf{266} vol. 4, (2014), 2590-2609.
\bibitem{BNSSW} J.H. Brown, G. Nagy, S. Reznikoff, A. Sims and D.P. Williams, \emph {Cartan subalgebras in C*-algebras of Hausdorff etale groupoids}, (2015), arXiv:1503.03521v2 
\bibitem{FGKP} C. Farsi, E. Gillaspy, S. Kang and J. Packer, \emph{Wavelets and graph C*-algebras}, (2016) arXiv:1601.00061v1.
\bibitem{GLR} D. Gon{\c{c}}alves, H. Li, and D. Royer, \emph{Faithful representations of graph algebras via branching systems}, Canad. Math. Bull. \textbf{59} (2016), 95-103.
\bibitem{GR4} D. Gon{\c{c}}alves and D. Royer, \emph{Branching systems and representations of {C}ohn--{L}eavitt path algebras of separated graphs}, J. Algebra \textbf{422} (2015), 413--426.
\bibitem{MR2903145} D. Gon{\c{c}}alves and D. Royer, \emph{Graph {${\rm C}^*$}-algebras, branching systems and the {P}erron-{F}robenius operator}, J. Math. Anal. Appl. \textbf{391} (2012), 457--465.
\bibitem{GR3} D. Gon{\c{c}}alves and D. Royer, \emph{On the representations of {L}eavitt path algebras}, J. Algebra \textbf{333} (2011), 258--272.
\bibitem{GR} D. Gon{\c{c}}alves and D. Royer, \emph{Perron-{F}robenius operators and representations of the {C}untz-{K}rieger algebras for infinite matrices}, J. Math. Anal. Appl. \textbf{351} (2009), 811--818.
\bibitem{GRUltragraph} D. Gon{\c{c}}alves and D. Royer, \emph{Ultragraphs and shifts spaces over infinite alphabets}, (2015), arXiv:1510.04649v2.
\bibitem{MR2848777} D. Gon{\c{c}}alves and D. Royer, \emph{Unitary equivalence of representations of algebras associated with graphs, and branching systems}, Functional Analysis and its Applications \textbf{45} (2011), 45--59.
\bibitem{Katsura:CJM03} T. Katsura, \emph{The ideal structures of crossed products of {C}untz algebras by quasi-free actions of abelian groups}, Canad. J. Math. \textbf{55} (2003), 1302--1338.
\bibitem{KMST} T. Katsura, P.S. Muhly, A. Sims and M. Tomforde, \emph{Graph algebras, Exel-{L}aca algebras, and ultragraph algebras coincide up to Morita equivalence}, J. Fur Die Reine Ang. Math. \textbf{640} (2010), 135-165.
\bibitem{Kawamura} K. Kawamura, \emph{Polynomial representations of the Cuntz algebras arising from permutations I}, General theory, Lett. Math. Phys. \textbf{71} (2005), 149--158.
\bibitem{KumjianPaskEtAl:PJM98} A. Kumjian, D. Pask, and I. Raeburn, \emph{Cuntz-{K}rieger algebras of directed graphs}, Pacific J. Math. \textbf{184} (1998), 161--174.
\bibitem{LY}  A. Lasota and J.A. Yorke, \emph{Exact dynamical systems and the Frobenius-Perron operator}, Trans.Amer.Math.Soc. \textbf{273} (1982), 375--384.
\bibitem{Lawson} M.V. Lawson, \emph{Primitive partial permutation representations of the polycyclic monoids and branching funtion systems}, Period. Math. Hungar. \textbf{58} (2009), 189--207.
\bibitem{Lin:introductiontoclassification01} H. Lin, An introduction to the classification of amenable {$C\sp *$}-algebras, World Scientific Publishing Co. Inc., River Edge, NJ, 2001, xii+320.
\bibitem{Raeburn:Graphalgebras05} I. Raeburn, Graph algebras, Published for the Conference Board of the Mathematical Sciences, Washington, DC, 2005, vi+113.
\bibitem{RaeburnWilliams:Moritaequivalenceand98} I. Raeburn  and D.P. Williams, Morita equivalence and continuous-trace {$C\sp *$}-algebras, American Mathematical Society, Providence, RI, 1998, xiv+327.
\bibitem{Tomforde:JOT03} M. Tomforde, \emph{A unified approach to {E}xel-{L}aca algebras and {$C\sp \ast$}-algebras associated to graphs}, J. Operator Theory \textbf{50} (2003), 345--368.
\bibitem{Tomforde:IUMJ03} M. Tomforde, \emph{Simplicity of ultragraph algebras}, Indiana Univ. Math. J. \textbf{52} (2003), 901--925.
\bibitem{Szyma'nski:IJM02} W. Szyma{\'n}ski, \emph{General {C}untz-{K}rieger uniqueness theorem}, Internat. J. Math. \textbf{13} (2002), 549--555.
\end{thebibliography}
\end{document}